\documentclass[reqno,11pt]{amsart}
\usepackage{amsmath}
\usepackage{amsxtra}
\usepackage{amscd}
\usepackage{amsthm}
\usepackage{amsfonts}
\usepackage{amssymb}
\usepackage{eucal}

\theoremstyle{plain}
\newtheorem{Thm}[subsection]{Theorem}
\newtheorem{Cor}[subsection]{Corollary}
\newtheorem{Lem}[subsection]{Lemma}
\newtheorem{Prop}[subsection]{Proposition}
\newtheorem{Conj}[subsection]{Conjecture}

\theoremstyle{definition}
\newtheorem{Def}[subsection]{Definition}

\theoremstyle{remark}

\newtheorem{Rem}[subsection]{Remark}

\newtheorem{conj}{Conjecture}

\numberwithin{equation}{section}
\renewcommand{\rm}{\normalshape}

\newif\ifShowLabels
\ShowLabelstrue
\newdimen\theight
\def\TeXref#1{%
    \leavevmode\vadjust{\setbox0=\hbox{{\tt
        \quad\quad  {\small \rm #1}}}%
    \theight=\ht0
    \advance\theight by \lineskip
    \kern -\theight \vbox to
    \theight{\rightline{\rlap{\box0}}%
    \vss}%
    }}%

\ShowLabelsfalse

\renewcommand{\sec}[2]{\section{#2}\label{S:#1}%
    \ifShowLabels \TeXref{{S:#1}} \fi}
\newcommand{\ssec}[2]{\subsection{#2}\label{SS:#1}%
    \ifShowLabels \TeXref{{SS:#1}} \fi}

\newcommand{\sssec}[2]{\subsubsection{#2}\label{SSS:#1}%
    \ifShowLabels \TeXref{{SSS:#1}} \fi}

\newcommand{\refs}[1]{Section ~\ref{S:#1}}
\newcommand{\refss}[1]{Section ~\ref{SS:#1}}
\newcommand{\refsss}[1]{Section ~\ref{SSS:#1}}
\newcommand{\reft}[1]{Theorem ~\ref{T:#1}}
\newcommand{\refl}[1]{Lemma ~\ref{L:#1}}
\newcommand{\refp}[1]{Proposition ~\ref{P:#1}}

\newcommand{\refe}[1]{\eqref{E:#1}}
\newcommand{\refco}[1]{Conjecture ~\ref{Co:#1}}

\newenvironment{thm}[1]%
    { \begin{Thm} \label{T:#1}  \ifShowLabels \TeXref{T:#1} \fi }%
    { \end{Thm} }

\renewcommand{\th}[1]{\begin{thm}{#1} \sl }
\renewcommand{\eth}{\end{thm} }

\newenvironment{lemma}[1]%
    { \begin{Lem} \label{L:#1}  \ifShowLabels \TeXref{L:#1} \fi }%
    { \end{Lem} }
\newcommand{\lem}[1]{\begin{lemma}{#1} \sl}
\newcommand{\elem}{\end{lemma}}

\newenvironment{propos}[1]%
    { \begin{Prop} \label{P:#1}  \ifShowLabels \TeXref{P:#1} \fi }%
    { \end{Prop} }
\newcommand{\prop}[1]{\begin{propos}{#1}\sl }
\newcommand{\eprop}{\end{propos}}

\newenvironment{corol}[1]%
    { \begin{Cor} \label{C:#1}  \ifShowLabels \TeXref{C:#1} \fi }%
    { \end{Cor} }
\newcommand{\cor}[1]{\begin{corol}{#1} \sl }
\newcommand{\ecor}{\end{corol}}

\newenvironment{defeni}[1]%
    { \begin{Def} \label{D:#1}  \ifShowLabels \TeXref{D:#1} \fi }%
    { \end{Def} }
\newcommand{\defe}[1]{\begin{defeni}{#1} \sl }
\newcommand{\edefe}{\end{defeni}}

\newenvironment{remark}[1]%
    { \begin{Rem} \label{R:#1}  \ifShowLabels \TeXref{R:#1} \fi }%
    { \end{Rem} }
\newcommand{\rem}[1]{\begin{remark}{#1}}
\newcommand{\erem}{\end{remark}}

\newenvironment{conjec}[1]%
    { \begin{Conj} \label{Co:#1}  \ifShowLabels \TeXref{Co:#1} \fi }%
    { \end{Conj} }
\renewcommand{\conj}[1]{\begin{conjec}{#1} \sl }
\newcommand{\econj}{\end{conjec}}

\newcommand{\eq}[1]%
    { \ifShowLabels \TeXref{E:#1} \fi
       \begin{equation} \label{E:#1} }
\newcommand{\eeq}{ \end{equation} }

\newcommand{\prf}{ \begin{proof} }
\newcommand{\epr}{ \end{proof} }

\newcommand\alp{\alpha}     

     \newcommand\Gam{\Gamma}

\newcommand\kap{\kappa}
\newcommand\lam{\lambda}        \newcommand\Lam{\Lambda}

\newcommand\calF{{\mathcal{F}}}

\newcommand\calK{{\mathcal{K}}}

\newcommand\calO{{\mathcal{O}}}

\newcommand\calR{{\mathcal{R}}}
\newcommand\calS{{\mathcal{S}}}
\newcommand\calT{{\mathcal{T}}}
\newcommand\calU{{\mathcal{U}}}

\newcommand\calW{{\mathcal{W}}}

\newcommand\PP{\mathbb{P}}
\renewcommand\AA{\mathbb{A}}

\newcommand\GG{\mathbb{G}}

\newcommand\ZZ{\mathbb{Z}}

\newcommand\CC{\mathbb{C}}

 \newcommand\grg{{\mathfrak{g}}}

 \newcommand\grl{{\mathfrak{l}}}

\newcommand\sdp{\times \hskip -0.3em {\raise 0.3ex
\hbox{$\scriptscriptstyle |$}}} 

\newcommand\Gr{\operatorname{Gr}}

\newcommand\Hom{\operatorname {Hom}}

\newcommand\Int{\operatorname{Int}}

\newcommand\Perv{\operatorname{Perv}}

\newcommand\SL{{\rm SL}}

\newcommand\olam{{\overline{\lambda}}}

\newcommand\hatG{{\widehat{G}}}

\newcommand\tilG{{\widetilde{G}}}

\newcommand\tilp{{\widetilde{p}}}

\newcommand\tilS{{\widetilde{S}}}

\newcommand\x{\times}
\newcommand\ten{\otimes}

\newcommand{\ra}{\rangle}
\newcommand{\la}{\langle}

\newcommand\nc{\newcommand}

\newcommand{\go}{G(\calO)}

\newcommand{\rc}{\rho^{\vee}}
\newcommand{\ogg}{{\overline \gg}}
\newcommand{\ogl}{\ogg^{\lam}}
\newcommand{\pgg}{\text{Perv}_{\go}(\gg)}
\newcommand{\IC}{{\operatorname{IC}}}

\newcommand{\iso}{{\widetilde\longrightarrow}}

\nc\aff{\operatorname{aff}}
\nc\oGr{\overline{\Gr}}
\nc\Bun{\operatorname{Bun}}
\nc\hgrg{\widehat{\grg}}
\renewcommand\Int{\operatorname{Int}}
\nc\bInt{\overline{\Int}}
\nc\hatLam{\widehat{\Lam}}
\nc\bmu{\overline{\mu}}
\nc\bnu{\overline{\nu}}
\nc\blambda{\overline{\lam}}
\renewcommand\SL{\operatorname{SL}}
\nc\ocalW{\overline{\calW}}
\nc\pos{\operatorname{pos}}
\nc\IH{\operatorname{IH}}
\nc\Rep{\operatorname{Rep}}
\nc\Gal{\operatorname{Gal}}
\nc{\tilGr}{\widetilde{\Gr}}

\nc\Pic{\operatorname{Pic}}

\nc{\HC}{{\mathcal{HC}}}
\nc{\on}{\operatorname}
\nc{\BA}{{\mathbb{A}}}
\nc{\BC}{{\mathbb{C}}}
\nc{\BM}{{\mathbb{M}}}
\nc{\BN}{{\mathbb{N}}}
\nc{\BP}{{\mathbb{P}}}
\nc{\BR}{{\mathbb{R}}}
\nc{\BZ}{{\mathbb{Z}}}
\nc{\BS}{{\mathbb{S}}}

\nc{\CA}{{\mathcal{A}}}
\nc{\CB}{{\mathcal{B}}}
\nc{\CalD}{{\mathcal D}}
\nc{\CE}{{\mathcal{E}}}
\nc{\CF}{{\mathcal{F}}}
\nc{\CG}{{\mathcal{G}}}
\nc{\CL}{{\mathcal{L}}}
\nc{\CM}{{\mathcal{M}}}
\nc{\CMM}{{\mathcal{M}^{\operatorname{gen}}_\hbar(-\rho)}}
\nc{\CN}{{\mathcal{N}}}
\nc{\CO}{{\mathcal{O}}}
\nc{\CP}{{\mathcal{P}}}
\nc{\CQ}{{\mathcal{Q}}}
\nc{\CR}{{\mathcal{R}}}
\nc{\CS}{{\mathcal{S}}}
\nc{\CT}{{\mathcal{T}}}
\nc{\CU}{{\mathcal{U}}}
\nc{\CV}{{\mathcal{V}}}
\nc{\CW}{{\mathcal{W}}}
\nc{\CX}{{\mathcal{X}}}
\nc{\CZ}{{\mathcal{Z}}}

\nc{\gen}{{\operatorname{gen}}}
\nc{\cM}{{\check{\mathcal M}}{}}
\nc{\csM}{{\check{\mathcal A}}{}}
\nc{\obM}{{\overset{\circ}{\mathbf M}}{}}
\nc{\oCA}{{\overset{\circ}{\mathcal A}}{}}
\nc{\obA}{{\overset{\circ}{\mathbf A}}{}}
\nc{\ooM}{{\overset{\circ}{M}}{}}
\nc{\osM}{{\overset{\circ}{\mathsf M}}{}}
\nc{\vM}{{\overset{\bullet}{\mathcal M}}{}}
\nc{\nM}{{\underset{\bullet}{\mathcal M}}{}}
\nc{\obD}{{\overset{\circ}{\mathbf D}}{}}
\nc{\cp}{{\overset{\circ}{\mathbf p}}{}}
\nc{\ofZ}{{\overset{\circ}{\mathfrak Z}}{}}

\nc{\fa}{{\mathfrak{a}}}
\nc{\fb}{{\mathfrak{b}}}
\nc{\fg}{{\mathfrak{g}}}
\nc{\fgl}{{\mathfrak{gl}}}
\nc{\fh}{{\mathfrak{h}}}
\nc{\fj}{{\mathfrak{j}}}
\nc{\fm}{{\mathfrak{m}}}
\nc{\fn}{{\mathfrak{n}}}
\nc{\fu}{{\mathfrak{u}}}
\nc{\fp}{{\mathfrak{p}}}
\nc{\frr}{{\mathfrak{r}}}
\nc{\fs}{{\mathfrak{s}}}
\nc{\fT}{{\mathfrak{T}}}
\nc{\ofT}{{\overline{\mathfrak T}}}
\nc{\ofS}{{\overline{\mathfrak S}}}
\nc{\fsl}{{\mathfrak{sl}}}
\nc{\hsl}{{\widehat{\mathfrak{sl}}}}
\nc{\hgl}{{\widehat{\mathfrak{gl}}}}
\nc{\hg}{{\widehat{\mathfrak{g}}}}
\nc{\chg}{{\widehat{\mathfrak{g}}}{}^\vee}
\nc{\hn}{{\widehat{\mathfrak{n}}}}
\nc{\chn}{{\widehat{\mathfrak{n}}}{}^\vee}

\nc{\fA}{{\mathfrak{A}}}
\nc{\fB}{{\mathfrak{B}}}
\nc{\fD}{{\mathfrak{D}}}
\nc{\fE}{{\mathfrak{E}}}
\nc{\fF}{{\mathfrak{F}}}
\nc{\fG}{{\mathfrak{G}}}
\nc{\fK}{{\mathfrak{K}}}
\nc{\fL}{{\mathfrak{L}}}
\nc{\fM}{{\mathfrak{M}}}
\nc{\fN}{{\mathfrak{N}}}
\nc{\frP}{{\mathfrak{P}}}
\nc{\fS}{{\mathfrak S}}
\nc{\fU}{{\mathfrak{U}}}
\nc{\fZ}{{\mathfrak{Z}}}

\nc{\bb}{{\mathbf{b}}}
\nc{\bc}{{\mathbf{c}}}
\nc{\be}{{\mathbf{e}}}
\nc{\bj}{{\mathbf{j}}}
\nc{\bn}{{\mathbf{n}}}
\nc{\bp}{{\mathbf{p}}}
\nc{\bq}{{\mathbf{q}}}
\nc{\bv}{{\mathbf{v}}}
\nc{\bx}{{\mathbf{x}}}
\nc{\by}{{\mathbf{y}}}
\nc{\bw}{{\mathbf{w}}}
\nc{\bA}{{\mathbf{A}}}
\nc{\bB}{{\mathbf{B}}}
\nc{\bC}{{\mathbf{C}}}
\nc{\bK}{{\mathbf{K}}}
\nc{\bD}{{\mathbf{D}}}
\nc{\bH}{{\mathbf{H}}}
\nc{\bM}{{\mathbf{M}}}
\nc{\bN}{{\mathbf{N}}}
\nc{\bT}{{\mathbf{T}}}
\nc{\bV}{{\mathbf{V}}}
\nc{\bW}{{\mathbf{W}}}
\nc{\bX}{{\mathbf{X}}}
\nc{\bP}{{\mathbf{P}}}
\nc{\bZ}{{\mathbf{Z}}}

\nc{\sA}{{\mathsf{A}}}
\nc{\sB}{{\mathsf{B}}}
\nc{\sC}{{\mathsf{C}}}
\nc{\sD}{{\mathsf{D}}}
\nc{\sF}{{\mathsf{F}}}
\nc{\sK}{{\mathsf{K}}}
\nc{\sM}{{\mathsf{M}}}
\nc{\sO}{{\mathsf{O}}}
\nc{\sQ}{{\mathsf{Q}}}
\nc{\sP}{{\mathsf{P}}}
\nc{\sZ}{{\mathsf{Z}}}
\nc{\sfp}{{\mathsf{p}}}
\nc{\sr}{{\mathsf{r}}}
\nc{\sfb}{{\mathsf{b}}}
\nc{\sfc}{{\mathsf{c}}}
\nc{\sd}{{\mathsf{d}}}
\nc{\sfl}{{\mathsf{l}}}

\nc{\BK}{{\bar{K}}}

\nc{\tA}{{\widetilde{\mathbf{A}}}}
\nc{\tB}{{\widetilde{\mathcal{B}}}}
\nc{\tg}{{\widetilde{\mathfrak{g}}}}
\nc{\tG}{{\widetilde{G}}}
\nc{\TM}{{\widetilde{\mathbb{M}}}{}}
\nc{\tO}{{\widetilde{\mathsf{O}}}{}}
\nc{\tU}{{\widetilde{\mathfrak{U}}}{}}
\nc{\TZ}{{\tilde{Z}}}
\nc{\tx}{{\tilde{x}}}
\nc{\tbv}{{\tilde{\bv}}}
\nc{\tfP}{{\widetilde{\mathfrak{P}}}{}}
\nc{\tz}{{\tilde{\zeta}}}
\nc{\tmu}{{\tilde{\mu}}}

\nc{\urho}{\underline{\rho}}
\nc{\uB}{\underline{B}}
\nc{\uC}{{\underline{\mathbb{C}}}}
\nc{\ui}{\underline{i}}
\nc{\ofP}{{\overline{\mathfrak{P}}}}

\nc{\hrho}{{\hat{\rho}}}
\nc{\Blambda}{{\boldsymbol{\lambda}}}
\nc{\Bblambda}{{\overline{\boldsymbol{\lambda}}}}

\nc{\unl}{\underline}
\nc{\ol}{\overline}
\nc{\one}{{\mathbf{1}}}
\nc{\two}{{\mathbf{t}}}

\nc{\Tot}{{\mathop{\operatorname{\rm Tot}}}}
\nc{\Hilb}{{\mathop{\operatorname{\rm Hilb}}}}
\nc{\CHom}{{\mathop{\operatorname{{\mathcal{H}}\it om}}}}
\nc{\defi}{{\mathop{\operatorname{\rm def}}}}
\nc{\length}{{\mathop{\operatorname{\rm length}}}}

\nc{\Cliff}{{\mathsf{Cliff}}}
\nc{\Fl}{{\mathsf{Fl}}}
\nc{\Fib}{{\mathsf{Fib}}}
\nc{\Coh}{{\mathsf{Coh}}}
\nc{\FCoh}{{\mathsf{FCoh}}}

\nc{\reg}{{\text{\rm reg}}}

\nc{\cplus}{{\mathbf{C}_+}}
\nc{\cminus}{{\mathbf{C}_-}}
\nc{\cthree}{{\mathbf{C}_*}}
\nc{\Qbar}{{\bar{Q}}}

\nc{\bh}{{\bar{h}}}
\nc{\bOmega}{{\overline{\Omega}}}
\nc\tGr{\widetilde{\Gr}}

\nc{\seq}[1]{\stackrel{#1}{\sim}}

\renewcommand\gg{\Gr_G}
\nc\uS{\underline{S}}
\begin{document}
\title{Pursuing the double affine Grassmannian II: Convolution}
\author{Alexander Braverman and Michael Finkelberg}

\begin{abstract}This is the second paper of a series (started by \cite{BF})
which describes
a conjectural analog of the affine Grassmannian for affine Kac-Moody groups
(also known as the double affine Grassmannian).
The current paper is dedicated to describing a conjectural analog of the
convolution diagram for the double affine Grassmannian.
In the case when $G=SL(n)$ our conjectures can be derived from \cite{Na}.
\end{abstract}
\maketitle
\sec{int}{Introduction}

\ssec{21}{The usual affine Grassmannian}Let $G$ be a connected complex reductive group and
let $\calK=\CC((s))$, $\calO=\CC[[s]]$. By the {\it affine
Grassmannian} of $G$ we shall mean the quotient
$\gg=G(\calK)/G(\calO)$. It is known (cf. \cite{BD,MV}) that
$\gg$ is the set of $\CC$-points of an ind-scheme over
$\CC$, which we will denote by the same symbol. Note that $\gg$ is defined
for any (not necessarily reductive) group $G$.

Let $\Lam=\Lam_G$ denote the coweight lattice of $G$ and let $\Lam^{\vee}$
denote the dual lattice (this is the weight lattice of $G$).
We let $2\rho_G^{\vee}$ denote the sum of the positive roots of $G$.

The group-scheme $G(\calO)$ acts on $\gg$ on the left and
its orbits can be described as follows.
One can identify the lattice $\Lam_G$ with
the quotient $T(\calK)/T(\calO)$. Fix $\lam\in\Lam_G$ and
let $s^{\lam}$ denote any lift of $\lam$ to $T(\calK)$.
Let $\gg^{\lam}$ denote the $\go$-orbit of $s^{\lam}$
(this is clearly independent of the choice of $s^{\lam}$).
The following result is well-known:
\lem{gras-orbits}
\begin{enumerate}
\item
$$
\gg=\bigcup\limits_{\lam\in\Lam_G}\gg^{\lam}.
$$
\item
We have $\Gr_G^{\lam}=\Gr_G^{\mu}$ if an only if $\lam$ and $\mu$ belong
to the same $W$-orbit on $\Lam_G$ (here $W$ is the Weyl group of $G$). In particular,
$$
\gg=\bigsqcup\limits_{\lam\in\Lam^+_G}\gg^{\lam}.
$$
\item
For every $\lam\in\Lam^+$ the orbit
$\gg^{\lam}$ is finite-dimensional and its dimension is
equal to $\la\lam,2\rho_G^{\vee}\ra$.
\end{enumerate}
\elem
Let $\ogl$ denote the closure of $\gg^{\lam}$ in $\gg$;
this is an irreducible projective algebraic variety; one has
$\gg^{\mu}\subset \ogl$ if and only if $\lam-\mu$ is a sum of positive roots of
$G^{\vee}$.
We will denote by $\IC^{\lambda}$ the intersection
cohomology complex on $\ogl$. Let $\pgg$ denote the category of
$G(\calO)$-equivariant perverse sheaves on $\gg$. It is known
that every object of this category is a direct sum of the
$\IC^{\lam}$'s.


\ssec{trans-finite}{Transversal slices}Consider the group $G[s^{-1}]\subset G((s))$; let us denote by
$G[s^{-1}]_1$ the kernel of the natural (``evaluation at $\infty$") homomorphism
$G[s^{-1}]\to G$. For any $\lam\in\Lam$ let $\Gr_{G,\lam}=G[s^{-1}]\cdot s^{\lam}$. Then
it is easy to see that one has
$$
\gg=\bigsqcup\limits_{\lam\in\Lam^+}\Gr_{G,\lam}
$$

Let also $\calW_{G,\lam}$ denote the $G[s^{-1}]_1$-orbit of $s^{\lam}$.
For any $\lam,\mu\in\Lam^+$, $\lam\geq \mu$ set
$$
\Gr^{\lam}_{G,\mu}=\gg^{\lam}\cap \Gr_{G,\mu},\quad
\oGr^{\lam}_{G,\mu}=\oGr_G^{\lam}\cap \Gr_{G,\mu}
$$
and
$$
\calW^{\lam}_{G,\mu}=\gg^{\lam}\cap \calW_{G,\mu},\quad
\ocalW^{\lam}_{G,\mu}=\oGr_G^{\lam}\cap \calW_{G,\mu}.
$$
Note that $\ocalW^{\lam}_{G,\mu}$ contains the point $s^{\mu}$ in it. The variety $\ocalW^{\lam}_{G,\mu}$ can be thought of
as a transversal slice to $\gg^{\mu}$ inside $\oGr_G^{\lam}$ at the point $s^{\mu}$ (cf. \cite{BF}, Lemma 2.9).
\ssec{convolution}{The convolution}
We can regard $G(\calK)$ as a total space of a $G(\calO)$-torsor over $\gg$. In particular,
by viewing another copy of $\gg$ as a $G(\calO)$-scheme, we can form the associated fibration
$$
\gg\star\gg:=G(\calK)\underset{\go}\x\gg.
$$

One has the natural maps $p,m:\gg\star\gg\to\gg$ defined as follows. Let $g\in G(\calK),x\in\gg$. Then
$$
p(g\times x)=g\, \text{mod}\,\go; \quad m(g\times x)=g\cdot x.
$$

For any $\lam_1,\lam_2\in\Lam_G^+$ let us set $\gg^{\lam_1}\star\gg^{\lam_2}$ to be the corresponding subscheme
of $\gg\star\gg$; this is a fibration over $\gg^{\lam_1}$ with the typical fiber $\gg^{\lam_2}$. Its closure is
$\oGr^{\lam_1}\star\oGr^{\lam_2}$. In addition, we define
$$
(\gg^{\lam_1}\star\gg^{\lam_2})^{\lam_3}=m^{-1}(\gg^{\lam_3})\cap
(\gg^{\lam_1}\star\gg^{\lam_2}).
$$
It is known (cf. \cite{Lu-qan}) that
\eq{Lusestimate}
\dim((\gg^{\lam_1}\star\gg^{\lam_2})^{\lam_3})=
\la \lam_1+\lam_2+\lam_3,\rc_G\ra.
\end{equation}
(It is easy to see that although $\rc_G\in \frac{1}{2}\Lam_G^{\vee}$,
the RHS of \refe{Lusestimate} is an integer whenever the above intersection
is non-empty.)

Starting from any perverse sheaf $\calT$ on $\gg$ and a $G(\calO)$-equivariant perverse
sheaf $\calS$ on $\gg$, we can form their twisted external product $\calT\widetilde{\boxtimes}\calS$, which will
be a perverse sheaf on $\gg\star\gg$. For two objects $\calS_1,\calS_2\in\pgg$ we define their convolution
$$
\calS_1\star\calS_2=m_!(\calS_1\widetilde{\boxtimes}\calS_2).
$$

The following theorem, which is a categorical version of the Satake equivalence,
is a starting point for this paper, cf. \cite{Lu-qan},\cite{Gi} and \cite{MV}. The best reference
so far is \cite{BD}, Sect. 5.3.

\th{grassmannian}
\begin{enumerate}
\item
Let $\calS_1,\calS_2\in \Perv_{\go}(\gg)$. Then
$\calS_1\star\calS_2\in \Perv_{\go}(\gg)$.
\item
The convolution $\star$ extends to a structure of
a tensor category on $\Perv_{\go}(\gg)$.
\item
As a tensor category, $\Perv_{G(\calO)}(\Gr_G)$ is
equivalent to the category $\text{Rep}(G^{\vee})$. Under this
equivalence, the object $\IC^{\lam}$ goes over to the irreducible representation
$L(\lam)$ of $G^{\vee}$ with highest weight $\lam$).
\end{enumerate}
\eth

\ssec{}{$n$-fold convolution}Similarly to the above, we can define the $n$-fold
convolution diagram
$$
m_n:\underset{n}{\underbrace{\Gr_G\star\cdots\star\Gr_G}}\to \Gr_G.
$$
Here
$$
\underset{n}{\underbrace{\Gr_G\star\cdots\star\Gr_G}}=\underset{n-1}{\underbrace{
G(\calK)\underset{G(\calO)}\x\cdots \underset{G(\calO)}\x G(\calK)}}\underset{G(\calO)}\x \Gr_G
$$
and $m_n$ is the multiplication map. Thus, given $n$ objects $\calS_1,\cdots,\calS_n$ of
$\Perv_{G(\calO)}(\Gr_G)$ we may consider the convolution $\calS_1\star\cdots\star\calS_n$; this
will be again an object of $\Perv_{G(\calO)}(\Gr_G)$ which under the equivalence of
\reft{grassmannian} corresponds to $n$-fold tensor product in $\Rep(G^{\vee})$. In particular,
let $\lam_1,\cdots,\lam_n$ be elements of $\Lam^+$. One can consider the corresponding
subvariety $\oGr_G^{\lam_1}\star\cdots\star\oGr_G^{\lam_n}$ in $\underset{n}{\underbrace{\Gr_G\star\cdots\star\Gr_G}}$.
Then the convolution $\IC^{\lam_1}\star\cdots\star\IC^{\lam_n}$ is just
the direct image $(m_n)_!(\IC(\oGr_G^{\lam_1}\star\cdots\star\oGr_G^{\lam_n}))$. In particular, we have
an isomorphism
\eq{decomp-transversal}
(m_n)_!(\IC(\oGr_G^{\lam_1}\star\cdots\star\oGr_G^{\lam_n}))\simeq
\bigoplus\limits_{\nu\in\Lam^+} \IC^{\nu}\ten \Hom(L(\nu),L(\lam_1)\ten\cdots\ten L(\lam_n)).
\end{equation}

\ssec{gaff}{The group $G_{\aff}$}From now on we assume that $G$ is almost simple and simply connected.
To a connected reductive group $G$ as above one can associate the corresponding
affine Kac-Moody group $G_{\aff}$ in the following way.
One can consider the polynomial loop group $G[t,t^{-1}]$ (this is an infinite-dimensional group ind-scheme)

It is well-known that $G[t,t^{-1}]$ possesses a canonical central extension $\tilG$ of $G[t,t^{-1}]$:
$$
1\to \GG_m\to \tilG\to G[t,t^{-1}]\to 1.
$$
Moreover,\ $\tilG$ has again a natural structure of a group ind-scheme.

The multiplicative group $\GG_m$ acts naturally on $G[t,t^{-1}]$ and this action lifts to $\tilG$.
We denote the corresponding semi-direct product by $G_{\aff}$; we also let $\grg_{\aff}$ denote its Lie algebra.

The Lie algebra $\grg_{\aff}$ is an untwisted affine Kac-Moody Lie algebra.
In particular,
it can be described by the corresponding affine root system. We denote by $\grg_{\aff}^{\vee}$ the
{\em Langlands dual affine Lie algebra} (which corresponds to the dual affine root system)
and by $G^{\vee}_{\aff}$ the corresponding dual affine Kac-Moody group, normalized by the property
that it contains $G^{\vee}$ as a subgroup (cf. \cite{BF}, Subsection 3.1 for more details).

We denote by $\Lam_{\aff}=\ZZ\x\Lam\x\ZZ$ the coweight lattice of $G_{\aff}$; this is the same as the
weight lattice of $G_{\aff}^{\vee}$. Here the first $\ZZ$-factor
is responsible for the center of $G_{\aff}^{\vee}$ (or $\hatG^{\vee}$);
it can also be thought of as coming from the loop
rotation in $G_{\aff}$. The second $\ZZ$-factor is responsible for the loop rotation in $G_{\aff}^{\vee}$  it may also be thought of
as coming from the center of $G_{\aff}$).
We denote by $\Lam_{\aff}^+$ the set of dominant weights of $G_{\aff}^{\vee}$ (which is the same as the set of dominant
coweights of $G_{\aff}$). We also denote by $\Lam_{\aff,k}$ the set of weights of $G_{\aff}^{\vee}$ of level $k$,
i.e. all the weights of the form $(k,\olam,n)$. We put $\Lam_{\aff,k}^+=\Lam_{\aff}^+\cap \Lam_{\aff,k}$.

\smallskip
\noindent
{\bf Important notational convention:}
From now on we shall denote elements
of $\Lam$ by $\blambda,\bmu...$ (instead of just writing $\lam,\mu...$ in order to distinguish them from
the coweights of $G_{\aff}$ (= weights of $G_{\aff}^{\vee}$), which we shall just denote by
$\lam,\mu...$

Let
$\Lam_k^+\subset \Lam$ denote the set of dominant coweights of $G$ such that $\la \blambda,\alp)\leq k$
when $\alp$ is the highest root of $\grg$.
Then it is well-known that a weight $(k,\olam,n)$ of
$G_{\aff}^{\vee}$ lies in $\Lam_{\aff,k}^+$  if and only if $\olam\in\Lam_k^+$ (thus $\Lam_{\aff,k}=\Lam_k^+\x \ZZ$).

Let also $W_{\aff}$ denote affine Weyl group of $G$ which is the semi-direct product of $W$ and $\Lam$.
It acts on the lattice $\Lam_{\aff}$ (resp. $\hatLam$)
preserving each $\Lam_{\aff,k}$ (resp. each $\hatLam_k$). In order to describe this action explicitly it is convenient
to set $W_{\aff,k}=W\ltimes k\Lam$ which naturally acts on $\Lam$. Of course the groups $W_{\aff,k}$ are canonically
isomorphic to $W_{\aff}$ for all $k$. Then the restriction of the $W_{\aff}$-action to $\Lam_{\aff,k}\simeq\Lam\x\ZZ$
comes from the natural $W_{\aff,k}$-action on the first multiple.

It is well known that every $W_{\aff}$-orbit on $\Lam_{\aff,k}$ contains unique element of $\Lam_{\aff,k}^+$.
This is equivalent to saying that $\Lam_k^+\simeq \Lam/W_{\aff,k}$.
\ssec{}{Transversal slices for $\Gr_{G_{\aff}}$}
Our main dream is to create an analog of the affine Grassmannian $\Gr_G$ and the above results
about it in the case when $G$ is replaced by the (infinite-dimensional) group $G_{\aff}$.
The first attempt to do so was made in \cite{BF}: namely, in {\em loc. cit.} we have constructed analogs of the varieties $\ocalW^{\lam}_{G,\mu}$ in the case
when $G$ is replaced by  $G_{\aff}$. In the current paper, we are going to construct analogs of the varieties
$m_n^{-1}(\ocalW^{\lam}_{G,\mu})\cap (\Gr_G^{\lam_1}\star\cdots\star\Gr_G^{\lam_n})$ and
$m_n^{-1}(\ocalW^{\lam}_{G,\mu})\cap (\oGr_G^{\lam_1}\star\cdots\star\oGr_G^{\lam_n})$ (here
$\lam=\lam_1+\cdots+\lam_n$) when
$G$ is replaced by $G_{\aff}$. We shall also construct (cf. \refss{fam}) analogs of the corresponding
pieces in the {\em Beilinson-Drinfeld Grassmannian} for $G_{\aff}$ (cf. \refss{fam-finite} for
a short digression on the Beilinson-Drinfeld Grassmannian for $G$).

To formulate the idea of our construction, let us first recall the construction of the affine analogs of
the varieties $\ocalW^{\lam}_{G,\mu}$.
Let $\Bun_G(\AA^2)$ denote the moduli space of principal $G$-bundles on $\PP^2$ trivialized at the ``infinite"
line
$\PP^1_{\infty}\subset \PP^2$. This is an algebraic variety which has connected components parametrized by non-negative
integers, corresponding to different values of the second Chern class of the corresponding bundles; we denote the corresponding
connected component by $\Bun_G^a(\AA^2)$ (here $a\geq 0$).
According to \cite{BFG} one can embed $\Bun_G^a(\AA^2)$ (as an open dense subset) into a larger variety
$\calU^a_G(\AA^2)$ which is called {\em the Uhlenbeck moduli space of $G$-bundles on $\AA^2$ of second
Chern class $a$}.\footnote{This space is an algebraic analog of the Uhlenbeck
compactification of the moduli space of instantons
on a Riemannian 4-manifold.}
Furthermore, for any $k\geq 0$, let $\Gam_k\subset SL(2)$ be the group of $k$-th roots of unity.
This group acts naturally on $\AA^2$ and $\PP^2$ and this action can be lifted to an action of
$\Gam_k$ on $\Bun_G(\AA^2)$ and $\calU_G(\AA^2)$. This lift depends on a choice of a homomorphism
$\Gam_k\to G$ which is responsible for the action of $\Gam_k$ on the trivialization of our $G$-bundles on
$\PP^1_{\infty}$; it is explained in \cite{BF} that to such a homomorphism one can associate a dominant
weight $\mu$ of $G_{\aff}^{\vee}$ of level $k$; in the future we shall denote the set of all such weights
by $\Lam^+_k$. We denote by $\Bun_{G,\mu}(\AA^2/\Gam_k)$ the set of fixed points of
$\Gam_k$ on $\Bun_G(\AA^2)$. In \cite{BF} we construct a bijection between connected components
of $\Bun_{G,\mu}(\AA^2/\Gam_k)$ and dominant weights $\lam$ of $G_{\aff}$ such that $\lam\geq \mu$.
We denote the corresponding connected component by $\Bun^{\lam}_{G,\mu}(\AA^2/\Gam_k)$; we also denote by
$\calU^{\lam}_{G,\mu}(\AA^2/\Gam_k)$ its closure in
$\calU_G(\AA^2)$.\footnote{More precisely, in \cite{BF} we construct an open and closed subvariety $\Bun^{\lam}_{G,\mu}(\AA^2/\Gam_k)$ inside
$\Bun_G(\AA^2/\Gam_k)$ and formulate a conjecture, saying that it is connected (and thus it is a connected
component of $\Bun_G(\AA^2/\Gam_k)$). This conjecture is proved in \cite{BF} for $G=SL(n)$ and it is still
open in general.} In \cite{BF} we explain in what sense the variety $\calU^{\lam}_{G,\mu}$ should be thought
of as the correct version of $\ocalW^{\lam}_{G_{\aff},\mu}$.
\ssec{BMKS}{Bundles on mixed Kleinian stacks}
For a scheme $X$ endowed with an action of a finite group $\Gam$ we shall denote by $X/\!/\Gam$ the scheme-theoretic (categorical)
quotient of $X$ by $\Gam$; similarly, we denote by $X/\Gam$ the corresponding quotient stack.\footnote{The categorical quotient $X/\!/\Gam$ may not exist in
general, but we will only deal with the case $X$ is affine or projective, when
the categorical quotient does exist.}

Given positive integers $k_1,\ldots,k_n$ such that  $\sum\limits_{i=1}^nk_i=k$,
we set
$$
\BA^2/\!/\Gamma_k=\unl{S}{}_k\subset\unl{\ol{S}}{}_k:=\BP^2/\!/\Gamma_k,\
\BA^2/\Gamma_k=S_k\subset\ol{S}_k=\BP^2/\Gamma_k.
$$
We define $\widetilde{S}_k$ (resp. $\widetilde{\ol{S}}_k$) as
the minimal resolution of $\unl{S}{}_k$ (resp. $\unl{\ol{S}}{}_k$) at the
point 0. The exceptional divisor $E\subset\widetilde{S}_k$ is an
$A_{k-1}$-diagram of projective lines $E_1,\ldots,E_{k-1}$. Since any $E_i$
is a $-2$-curve, it is possible to blow down an arbitrary subset of
$\{E_1,\ldots,E_{k-1}\}$. We set ${\vec k}=(k_1,\ldots,k_n)$, and we
define $\unl{S}_{\vec k}$
(resp. $\unl{\ol{S}}{}_{\vec k}$) as the result of blowing down
all the lines except for $E_{k_1},E_{k_1+k_2},\ldots,E_{k_1+\ldots+k_{n-1}}$
in $\widetilde{S}_k$ (resp. $\widetilde{\ol{S}}{}_k$).
The surface $\unl{S}{}_{\vec k}$
(resp. $\unl{\ol{S}}{}_{\vec k}$) possesses canonical stacky resolution
$S_{\vec k}$ (resp. ${\ol{S}}{}_{\vec k}$). We will denote by
$s_1,\ldots,s_n\in S_{\vec k}$ the torus fixed points with the
automorphism groups $\Gamma_{k_1},\ldots,\Gamma_{k_n}$.

We denote by $\Bun_G(S_{\vec k})$ the moduli space of $G$-bundles
on $\ol{S}{}_{\vec k}$ trivialized on the boundary divisor
$\ol{S}{}_{\vec k}\setminus S_{\vec k}$.
For a bundle $\CF\in\Bun_G(S_{\vec k})$, the group
$\Gamma_{k_i}$ acts on its fiber $\CF_{s_{k_i}}$ at the point
$s_{k_i}$, and hence defines a conjugacy class of maps
$\Gamma_{k_i}\to G$, i.e. an element of $\Lambda_{k_i}^+$.
Similarly, the action of $\Gamma_k$ at the fiber of $\CF$ at
infinity defines an element of $\Lambda_k^+$. We denote by
$\Bun_{G,\bmu}^{\blambda{}^{(1)},\ldots,\blambda{}^{(n)}}(S_{\vec k})$
the subset of $\Bun_G(S_{\vec k})$ formed by all
$\CF\in\Bun_G(S_{\vec k})$ such that $\CF_{s_{k_i}}$ is of
class $\blambda{}^{(i)}$, and $\CF_\infty$ is of class $\bmu$.
To unburden the notations, we will write $\Bblambda$ for
$(\blambda{}^{(1)},\ldots,\blambda{}^{(n)})$, and
$\Bun_{G,\bmu}^\Bblambda(S_{\vec k})$ for
$\Bun_{G,\bmu}^{\blambda{}^{(1)},\ldots,\blambda{}^{(n)}}(S_{\vec k})$.
Clearly, it is a union of connected components of $\Bun_G(S_{\vec k})$.
We denote by
$\Bun_{G,\bmu}^{\Bblambda,d/k}(S_{\vec k})$
the intersection of
$\Bun_{G,\bmu}^\Bblambda(S_{\vec k})$
with $\Bun_G^{d/k}(S_{\vec k})$ ($G$-bundles of second Chern
class $d/k$. Here $d/k$ is the second Chern class on the stack; it is a
rational number with denominator $k$).
\ssec{uhl-conv}{Uhlenbeck spaces and convolution}
Our first goal in this paper is to define a certain partial Uhlenbeck
compactification $\CU_{G,\bmu}^{\Bblambda,d/k}
(S_{\vec k})
\supset\Bun_{G,\bmu}^{\Bblambda,d/k}(S_{\vec k})$.
The definition is given in \refs{SL} for
$G=\SL(N)$ using Nakajima's quiver varieties and in \refs{cotan} for general $G$ (using
all possible embeddings of $G$ into $\SL(N)$). Choosing certain lifts
$\lambda^{(i)}$ of $\blambda{}^{(i)}$ and $\mu$ of $\bmu$ to level $k$ dominant
weights of $G^\vee_{\operatorname{aff}}$ we will redenote
$\CU_{G,\bmu}^{\Bblambda,d/k}(S_{\vec k})$ by
$\CU_{G,\mu}^\Blambda(S_{\vec k})$. We will also
construct a proper birational morphism $\varpi:\
\CU_{G,\mu}^\Blambda(S_{\vec k})\to
\CU_{G,\mu}^\lambda(S_k)$ for $\lambda=\lambda^{(1)}+\ldots+\lambda^{(n)}$.
We believe that $\varpi$ is the correct analog of the convolution morphism
$$
m_n^{-1}(\ocalW^{\lam}_{G,\mu})\cap (\oGr_G^{\lam_1}\star\cdots\star\oGr_G^{\lam_n})\to
\ocalW^{\lam}_{G,\mu}.
$$
In particular, in the case $G=\SL(N)$ we prove an analog of \refe{decomp-transversal}
for the morphism $\varpi$ (the proof follows from the results of \cite{Na} by a fairly
easy combinatorial argument). We conjecture that a similar decomposition holds for general
$G$.

Let us note that the above conjecture is somewhat reminiscent of the results of \cite{BK} where similar
moduli spaces have been used in order to prove the existence of convolution for the {\em spherical
Hecke algebra} of $G_{\aff}$ (recall that the tensor category $\Perv_{G(\calO)}(\Gr_G)$ is a categorification of
the spherical Hecke algebra of $G$).
\ssec{axio}{Axiomatic approach to Uhlenbeck spaces}
We note that the above constructions of the relevant Uhlenbeck spaces and morphisms between them
are rather ad hoc; to give the reader certain
perspective, let us formulate what we would expect from Uhlenbeck spaces for general smooth
2-dimensional Deligne-Mumford stacks.
The constructions of \refs{SL} and \refs{cotan} may be viewed as a partial verification
of these expectations in the case of mixed Kleinian stacks.

\sssec{Spaces}{Spaces}Let $S$ be a smooth 2-dimensional Deligne-Mumford stack.
Let $\ol{S}$ be its smooth compactification, and let $D\subset\ol{S}$ be
the divisor at infinity. Let $\unl{S}\subset\unl{\ol{S}}$ be the coarse moduli
spaces. In our applications we only consider the stacks with cyclic
automorphism groups of points; more restrictively, only toric stacks.

Let $G$ as before be an almost simple simply connected complex algebraic group.
We assume that there are no $G$-bundles on $\ol{S}$ equipped with a
trivialization on $D$ with nontrivial automorphisms (preserving the
trivialization). In this case there is a fine moduli space $\Bun_G(S)$
of the pairs (a $G$-bundle on $\ol{S}$; its trivialization on $D$).
We believe that $\Bun_G(S)$ is open dense in the {\em Uhlenbeck completion}
$\CU_G(S)$. We believe that $\CU_{\SL(N)}(S)$ is a certain quotient of the
moduli stack of perverse coherent sheaves on $S$ which are $n$-dimensional
vector bundles off finitely many points.

A nontrivial homomorphism $\varrho:\ G\to\SL(N)$ gives rise to the closed
embedding $\varrho_*:\ \Bun_G(S)\hookrightarrow\Bun_{\SL(N)}(S)$ which we
expect to extend to a morphism $\CU_G(S)\hookrightarrow\CU_{\SL(N)}(S)$.
\sssec{morphisms}{Morphisms}
Assume that we have a proper morphism $\pi:\ \unl{\ol{S}}\to\unl{\ol{S}}{}'$ which is an
isomorphism in the neighbourhoods of $D,D'$. We believe $\pi$ gives rise
to a birational proper morphism $\varpi:\ \CU_G(S)\to\CU_G(S')$.
If $\varrho:\ G\to\SL(N)$ is a nontrivial representation of $G$, and
$\phi\in\CU_G(S)$, we choose a perverse coherent sheaf $F$ on $S$ representing
$\varrho_*(\phi)$. According to Theorem~4.2 of~\cite{Ka}, there is an
equivalence of derived coherent categories on $S$ and $S'$ (it is here that we
need the assumption that $S$ and $S'$ are toric).
This equivalence takes $F$ to a perverse coherent sheaf
$F'$ on $S'$. We believe that the class of $F'$ in $\CU_{\SL(N)}(S')$
equals $\varrho_*(\varpi(\phi))$.
\sssec{families}{Families}
Assume we have a morphism $\ol\CS\to\CX$ where $\CX$ is a variety, and for
every $x\in\CX$ the fiber $\ol{S}_x$ over $x$ is of type considered
in~\refsss{Spaces}. Then there should exist a morphism of varieties
$\CU_G(\CS)\to\CX$ such that for every $x\in\CX$ the fiber
$\CU_G(\CS)_x$ is isomorphic to $\CU_G(S_x)$ where $S_x\subset\ol{S}_x$ is the
canonical stacky resolution of $\unl{S}{}_x\subset\ol{\unl S}{}_x$,
cf.~\refss{resol-equiv} (note that we {\em do not require} the existence of a
family of stacks over $\CX$ with fibers $S_x$).

\ssec{}{Acknowledgements}The idea of using mixed Kleinian stacks in order to describe convolution in the
double affine Grassmannian was suggested to us by E.~Witten (a differential-geometric approach to this problem
is discussed in Section 5 of \cite{Wi}); we are very grateful to him for sharing his ideas with us and for
numerous very interesting conversations on the subject. Both authors would
like to thank H.~Nakajima for a lot of very helpful discussions and
in particular for his patient explanations of the contents
of~\cite{NaM},~\cite{Na}. We are obliged to R.~Bezrukavnikov,
D.~Kaledin, D.~Kazhdan, and A.~Kuznetsov for the useful discussions.
This paper was completed when the first author was visiting the J.-V.~Poncelet CNRS
laboratory at the Independent University of Moscow.

A.~B. was partially supported by the NSF grant DMS-0600851.
M.~F. was partially supported by the RFBR grant
09-01-00242, the HSE Science
Foundation award No. 11-09-0033,
the Ministry of Education and
Science of Russian Federation, grant No. 2010-1.3.1-111-017-029,
and the AG Laboratory HSE, RF government grant, ag. 11.G34.31.0023.

\sec{SL}{The case of $G=\SL(N)$}
\ssec{resol-equiv}{Stacky resolutions and derived equivalences}
In this Section we would like to implement the constructions announced
in~\refss{uhl-conv} in the case $G=\SL(N)$.
To do that let us first discuss some preparatory material.

Let $\uS$ be an algebraic surface and let $s_1,\ldots,s_n$ be distinct points on $\uS$ such that
the formal neighbourhood of $s_i$ is isomorphic to the formal neighbourhood of
$0$ in the surface $\AA^2/\!/\Gam_{k_i}$ for some $k_i\geq 1$; note that Artin's algebraization theorem
implies that such an isomorphism exists also \'etale-locally. Let us also assume that
$\uS$ is smooth away from $s_1,\ldots,s_n$. Recall that
for any $k\geq 1$ the surface $\AA^2/\!/\Gam_k$ possesses canonical minimal resolution
$\pi:\widetilde{\AA^2/\!/\Gam_k}\to \AA^2/\!/\Gam_k$ whose special fiber is a tree of type $A_{k-1}$ of $\PP^1$'s having
self-intersection $-2$. Similarly, we have a stacky resolution
$\AA^2/\Gam_k\to \AA^2/\!/\Gam_k$.
The existence of the above resolutions implies the existence of a resolution
$\tilS\to \uS$ and a stacky resolution\footnote{The existence follows from
the fact that every automorphism of $\BA^2/\Gam_k$ which is trivial over
$\AA^2/\!/\Gam_k\setminus\{0\}$, is trivial and the same is true over any
etale neighbourhood of $0$ in $\BA^2/\!/\Gam_k$.}
$S\to\uS$ which near every $s_i$ are \'etale locally isomorphic
to respectively $\widetilde{\AA^2/\!/\Gam_{k_i}}$ and $\AA^2/\Gam_{k_i}$.

For any scheme $Y$ let us denote by $D^b\Coh(Y)$ the bounded derived category of coherent sheaves on $Y$.
Recall (cf. \cite{KV} and \cite{BKR}) that we have an equivalence of derived categories
$$
\Psi:D^b\Coh(\widetilde{\AA^2/\!/\Gam_k})\to D^b\Coh(\AA^2/\Gam_k).
$$
This equivalence is given by a kernel which is a sheaf on $\widetilde{\AA^2/\!/\Gam_k}\x \AA^2/\Gam_k$ (and not
a complex of sheaves). Thus (by gluing in \'etale topology) a similar kernel can also be defined on the product
$\tilS\x S$ and it will define an equivalence
$D^b\Coh(\tilS)\to D^b\Coh(S)$ which we shall again denote by $\Psi$.
\ssec{zeta}{}
Recall the setup of~\cite{BF}~7.1.
Following~\cite{Na} we denote by
$I=\{1,\ldots,k\}$ the set of vertices of the affine cyclic quiver;
$k$ stands for the affine vertex, and $I_0=I\setminus\{k\}$.
Given ${\vec a}=(a_1,\ldots,a_n)\in\AA^n$ such that
\eq{ka}
k_1a_1+\ldots+k_na_n=0,
\end{equation}
we consider a
$k$-tuple $(b_1=a_1,\ldots,b_{k_1}=a_1,b_{k_1+1}=a_2,\ldots,b_{k_1+k_2}=a_2,
\ldots,b_{k_1+\ldots+k_{n-1}}=a_{n-1},b_{k_1+\ldots+k_{n-1}+1}=a_n,\ldots,
b_k=a_n)$. We consider another $k$-tuple of complex numbers
$\zeta^\circ_\BC$ such that $\zeta^\circ_{\BC,i}:=b_i-b_{i+1}$ for
$i=1,\ldots,k$ (for $i=k$ it is understood that $i+1=1$).

Furthermore, we set
$I_0\supset I_0^+:=\{k_1,k_1+k_2,\ldots,k_1+\ldots+k_{n-1}\}$,
and $I_0^0:=I_0\setminus I_0^+$. We consider a vector
$\zeta^\bullet_\BR\in\BR^I$
with coordinates $\zeta^\bullet_{\BR,i}=0$ for $i\in I_0^0$, and
$\zeta^\bullet_{\BR,j}=1$ for $j\in I_0^+$, and $\zeta^\bullet_{\BR,k}=1-n$.

Recall the setup
of~section~1 of~\cite{NaM}. In this note we are concerned with the cyclic
$A_{k-1}$-quiver only, so in particular, $\delta=(1,\ldots,1)$. We consider the
GIT quotient $X_{(\zeta^\circ_\BC,\zeta^\bullet_\BR)}:=\{\xi\in M(\delta,0):\
\mu(\xi)=-\zeta^\circ_\BC\}/\!/_{-\zeta^\bullet_\BR}(G_\delta/\BC^*)$
(see~(1.5) of~\cite{NaM}). It is a partial resolution of the categorical
quotient $X_{(\zeta^\circ_\BC,0)}:=\{\xi\in M(\delta,0):\
\mu(\xi)=-\zeta^\circ_\BC\}/\!/(G_\delta/\BC^*)$.
The above surfaces admit the following explicit description:
The surface $X_{(\zeta^\circ_\BC,0)}$ is isomorphic to the affine
surface given by the equation
$$
xy=(z-a_1)^{k_1}\ldots (z-a_n)^{k_n}.
$$
Note that when all $a_i$ are equal to $0$ we just get the equation $xy=z^n$ which defines
a surface isomorphic to $\AA^2/\!/\Gam_k$.
The surface $X_{(\zeta^\circ_\BC,\zeta^\bullet_\BR)}$ has the following properties:
if all the points $a_i$ are distinct, then
$X_{(\zeta^\circ_\BC,\zeta^\bullet_\BR)}=X_{(\zeta^\circ_\BC,0)}$.
On the other hand, if all $a_i$ are equal (and thus they have to be equal to zero by \refe{ka})
then $X_{(\zeta^\circ_\BC,\zeta^\bullet_\BR)}$ is obtained from $\widetilde{\AA^2/\!/\Gam_k}$ by blowing down
all the exceptional $\PP^1$'s except those whose numbers are $k_1,k_1+k_2,\ldots,k_1+\ldots+k_{n-1}$.
We leave the general case (i.e. the case of general $a_i$'s) to the reader.

The surface $X_{(\zeta^\circ_\BC,\zeta^\bullet_\BR)}$ (resp. $X_{(\zeta^\circ_\BC,0)}$) is of the type
discussed in \refss{resol-equiv} and thus it has canonical minimal stacky resolution, which we shall denote
by $S^{\vec a}_{\vec k}$ (resp.
$S'{}^{\vec a}_{\vec k}$).

If we choose a generic stability condition $\zeta^\circ_\BR$ in the hyperplane
$\zeta_\BR\cdot\delta=0$, then the corresponding GIT quotient
$X_{(\zeta^\circ_\BC,\zeta^\circ_\BR)}$ is smooth; moreover, it is the
minimal resolution of singularities of
$X_{(\zeta^\circ_\BC,\zeta^\bullet_\BR)}$. Recall the compactification
$\ol{X}_{(\zeta^\circ_\BC,\zeta^\circ_\BR)}$ introduced in~section~3
of~\cite{NaM}. According to \refss{resol-equiv} we have the equivalence
$\Psi:\ D^b\Coh(\ol{X}_{(\zeta^\circ_\BC,\zeta^\circ_\BR)})\iso
D^b\Coh(\ol{S}{}^{\vec a}_{\vec k})$.
Recall the line bundles $\CR_i,\ i\in I$, and their homomorphisms $\xi$ on
$\ol{X}_{(\zeta^\circ_\BC,\zeta^\circ_\BR)}$, introduced in sections~1(iii)
and~3(i) of~\cite{NaM}. We will denote $\Psi(\CR_i)$ by $\CR_i^\bullet$.
This is a line bundle on  $\ol{S}{}^{\vec a}_{\vec k}$
(this follows from the fact that a similar statement is true for
the equivalence $D^b\Coh(\widetilde{\AA^2/\!/\Gam_k})\to D^b\Coh(\AA^2/\Gam_k)$ under which the
bundle $\calR_i$ goes to the $\Gam_k$-equivariant sheaf on $\AA^2$ corresponding
to the structure sheaf of $\AA^2$ on which $\Gam_k$ acts by its $i$-th character).
\ssec{lambda}{} We consider the quiver variety
$\fM_{(\zeta^\circ_\BC,\zeta^\bullet_\BR)}(V,W)$ for the stability
condition $\zeta^\bullet_\BR$, see section~2 of~\cite{Na}.
We consider a vector
$\zeta^\pm_\BR:=\zeta^\bullet_\BR\pm(\varepsilon,\ldots,\varepsilon)
\in\BR^I$ for $0<\varepsilon\ll 1$. Note that it lies in an
(open) chamber of the stability conditions, so the corresponding
quiver varieties $\fM_{(\zeta^\circ_\BC,\zeta^\pm_\BR)}(V,W)$ are smooth.
Moreover, since
$\zeta^\bullet_\BR$ lies in a face adjacent to the chamber of
$\zeta^\pm_\BR$, we have the proper morphism
$\pi_{\zeta^\bullet_\BR,\zeta^\pm_\BR}:\
\fM_{(\zeta^\circ_\BC,\zeta^\pm_\BR)}(V,W)\to
\fM_{(\zeta^\circ_\BC,\zeta^\bullet_\BR)}(V,W)$.

The construction of~\cite{NaM}~(1.7),~3(ii) associates to any ADHM data
$(B,a,b)\in\bM(V,W)$ satisfying $\mu(B,a,b)=\zeta^\circ_\BC$
a complex of vector bundles
\eq{17}
\on{L}(\CR^{\bullet*},V)(-\ell_\infty)
\stackrel{\sigma}{\to}\on{E}(\CR^{\bullet*},V)
\oplus\on{L}(\CR^{\bullet*},W)
\stackrel{\tau}{\to}\on{L}(\CR^{\bullet*},V)(\ell_\infty)
\end{equation}
on $\ol{S}{}^{\vec a}_{\vec k}$.

The following proposition is a slight generalization of~Proposition~4.1
of~\cite{NaM}.

\prop{41}
Let $(B,a,b)\in\mu^{-1}(\zeta^\circ_\BC)$ and consider the complex~\refe{17}.
We consider $\sigma,\tau$ as linear maps on the fiber at a point in
$\ol{S}{}^{\vec a}_{\vec k}$. Then

(1) $(B,a,b)$ is $\zeta^-_\BR$-stable if and only if $\sigma$ is injective
possibly except finitely many points, and $\tau$ is surjective at any point.

(2) $(B,a,b)$ is $\zeta^\bullet_\BR$-semistable if and only if $\sigma$ is
injective and $\tau$ is surjective possibly except finitely many points.
\eprop

\prf The proof is parallel to that of~Proposition~4.1 of~\cite{NaM},
with the use of~Lemma~3.2 of~\cite{Na} in place of~Corollary~4.3
of~\cite{NaM}. \epr

\ssec{Levi}{}
We consider the Levi subalgebra $\grl\subset\fsl(k)\subset\fsl(k)_{\on{aff}}$
whose set of simple roots is
$I_0^0$, i.e. $\{\alpha_1,\ldots,\alpha_{k_1-1},\alpha_{k_1+1},\ldots,
\alpha_{k_1+k_2-1},\ldots,\alpha_{k-1}\}$. We will denote by $\BZ[I_0^0]$
the root lattice of $\grl$. The multiplication by the affine Cartan matrix
$A^{(1)}_{k-1}$ {\em embeds} $\BZ[I_0^0]$ into the weight lattice
$P_{\on{aff}}$ of $\fsl(k)_{\on{aff}}$ spanned by the fundamental weights
$\omega_0,\ldots,\omega_{k-1}$, so we will identify $\BZ[I_0^0]$ with a
sublattice of $P_{\on{aff}}$. The inclusion $\grl\subset\fsl(k)$ also gives
rise to the embedding $\BZ[I_0^0]\subset P$ into the weight lattice of
$\fsl(k)$.

We have ${\mathbf w}=\underline{\dim}W=(w_1,\ldots,w_k),\ {\mathbf
v}=\underline{\dim}V=(v_1,\ldots,v_k)$. We set $N:=w_1+\ldots+w_k$.
Recall the setup of~\cite{BF}~7.3. We associate to the pair
$(\bv,\bw)$ the $\fsl(k)_{\on{aff}}$-weight
$\bw'=\sum_{i=1}^kw'_i\omega_i:=
\sum_{i=1}^kw_i\omega_i-\sum_{i=1}^kv_i\alpha_i$.
In this note we restrict ourselves to the pairs $(\bv,\bw)$
satisfying the condition
\eq{uslovie}
\bw'\in N\omega_0+\BZ[I_0^0]
\end{equation}
The geometric meaning of this condition is as follows. The~\refp{41} implies
that $\fM_{(\zeta^\circ_\BC,\zeta^\bullet_\BR)}^{\on{reg}}(V,W)$ is the
moduli space of vector bundles on the stack
$\ol{S}{}^{\vec a}_{\vec k}$ trivialized at
infinity. The condition~\refe{uslovie} guaranties that these vector bundles
have trivial determinant, i.e. reduce to $\SL(N)$.

In effect, the determinant in question is a line bundle on
$\ol{S}{}^{\vec a}_{\vec k}$ trivialized at infinity.
So the determinant is trivial iff its restriction to the open substack
$S^{\vec a}_{\vec k}$ is trivial, i.e. is a zero element
of $\on{Pic}(S^{\vec a}_{\vec k})$.
Recall that $K(S^{\vec a}_{\vec k})\simeq P_{\on{aff}},\
\CR^\bullet_i\mapsto\omega_i$, and we have the homomorphism
$\det:\ K(S^{\vec a}_{\vec k})\to
\on{Pic}(S^{\vec a}_{\vec k})$. The class in
$K(S^{\vec a}_{\vec k})\simeq P_{\on{aff}}$ of any vector
bundle in $\fM_{(\zeta^\circ_\BC,\zeta^\bullet_\BR)}^{\on{reg}}(V,W)$ is
given by $\bw'\in P_{\on{aff}}$. So the triviality of its determinant is a
consequence of the following lemma.

\lem{Pic}
There is a canonical isomorphism
$\on{Pic}(S^{\vec a}_{\vec k})\simeq P/\BZ[I_0^0]$
such that the homomorphism $\det:\ K(S^{\vec a}_{\vec k})\to
\on{Pic}(S^{\vec a}_{\vec k})$ identifies with the
composition of the projection $P_{\on{aff}}\to P,\ \omega_i\mapsto\omega_i-
\delta_{i0}\omega_0$, and the projection $P\to P/\BZ[I_0^0]$.
\elem

\prf
Let $\widetilde{X}_{(\zeta^\circ_\BC,\zeta^\bullet_\BR)}$ stand for the
minimal resolution of the surface
$X_{(\zeta^\circ_\BC,\zeta^\bullet_\BR)}$.
Let $X^\circ_{(\zeta^\circ_\BC,\zeta^\bullet_\BR)}$ stand for the open
subset of $X_{(\zeta^\circ_\BC,\zeta^\bullet_\BR)}$ obtained by removing
all the singular points. The projection
$\widetilde{X}_{(\zeta^\circ_\BC,\zeta^\bullet_\BR)}\to
X_{(\zeta^\circ_\BC,\zeta^\bullet_\BR)}$ identifies
$X^\circ_{(\zeta^\circ_\BC,\zeta^\bullet_\BR)}$ with the open subset
of $\widetilde{X}_{(\zeta^\circ_\BC,\zeta^\bullet_\BR)}$ obtained by removing
the components $\{E_i,\ i\in I^0_0\}$ of the exceptional divisor.
Since any line bundle on $X^\circ_{(\zeta^\circ_\BC,\zeta^\bullet_\BR)}$
extends uniquely to a line bundle on $S^{\vec a}_{\vec k}$,
we obtain the restriction to the open subset homomorphism
$\on{Pic}(\widetilde{X}_{(\zeta^\circ_\BC,\zeta^\bullet_\BR)})
\twoheadrightarrow\on{Pic}(X^\circ_{(\zeta^\circ_\BC,\zeta^\bullet_\BR)})=
\on{Pic}(S^{\vec a}_{\vec k})$. Clearly, the kernel of this
restriction homomorphism is spanned by the classes of the line bundles
$\langle[\CO(E_i)],\ i\in I^0_0\rangle$ in
$\on{Pic}(\widetilde{X}_{(\zeta^\circ_\BC,\zeta^\bullet_\BR)})$.

Now recall that we have a canonical isomorphism
$\on{Pic}(\widetilde{X}_{(\zeta^\circ_\BC,\zeta^\bullet_\BR)})\simeq P$
such that the composition $\det:\ P_{\on{aff}}=
K(\widetilde{X}_{(\zeta^\circ_\BC,\zeta^\bullet_\BR)})\to
\on{Pic}(\widetilde{X}_{(\zeta^\circ_\BC,\zeta^\bullet_\BR)})\simeq P$
identifies with the projection
$p:\ P_{\on{aff}}\to P,\ \omega_i\mapsto\omega_i-
\delta_{i0}\omega_0$. Moreover, the class $[\CO(E_i)]\in
\on{Pic}(\widetilde{X}_{(\zeta^\circ_\BC,\zeta^\bullet_\BR)})$ gets identified
with $p(\alpha_i)$. This follows by embedding
$\widetilde{X}_{(\zeta^\circ_\BC,\zeta^\bullet_\BR)}$ as a slice into Grothendieck
simultaneous resolution $\widetilde{\mathfrak{sl}}_k$.

This completes the proof of the lemma.
\epr

\ssec{transpon}{}
Our next goal is to encode the quiver data $(\bv,\bw)$ by the
weight data of $\fsl(N)_{\on{aff}}$. From now on we assume that $\bw$
corresponds to an $N$-dimensional representation of $\Gamma_k$ with
{\em trivial determinant}, i.e. a homomorphism $\Gamma_k\to\SL(N)$.
Then the dominant weight $w_1\omega_1+\ldots+w_{k-1}\omega_{k-1}$ of
$\fsl(k)$ is actually a weight of $\on{PSL}(k)$, and can be written
uniquely as a generalized Young diagram $\tau=(\tau_1\geq\ldots\geq\tau_k)$
such that $\tau_i-\tau_{i+1}=w_i$ for any $1\leq i\leq k-1$, and
$\tau_1-\tau_k\leq N$, and $\tau_1+\ldots+\tau_k=0$, cf.~\cite{BF}~7.3.
Under the bijection $\Psi_{N,k}$ of {\em loc. cit.} $\bw$ corresponds to
a level $k$ dominant weight $\bmu\in\Lambda^+_k$ of $\widehat{\fsl(N)}$
which can also be written as a generalized Young diagram
$(\mu_1\geq\ldots\geq\mu_N)$ such that $\mu_1-\mu_N\leq k$, and
$\mu_1+\ldots+\mu_N=0$. We write $\tau={}^t\bmu$, and $\bmu={}^t\tau$.

Here is an explicit construction of the transposition operation on the
generalized Young diagrams. If $\bmu$ consists of all zeroes, then so does
$\tau$. Otherwise we assume $\mu_r>0\geq\mu_{r+1}$ for some $0<r<N$.
Then we have an {\em ordinary} Young diagram $\bmu':=
(k+\mu_{r+1}\geq k+\mu_{r+2}\geq
\ldots\geq k+\mu_N\geq\mu_1\geq\ldots\geq\mu_r)$ formed by {\em positive}
integers. We denote the {\em ordinary} transposition $^t\bmu'$ by $\tau'=
(\tau'_1\geq\ldots\geq\tau'_k)$, and finally we set $\tau={}^t\bmu:=
(\tau'_1+r-N\geq\ldots\geq\tau'_k+r-N)$. In other words,
\eq{javno}
\tau={}^t\bmu=(r^{\mu_r},
(r-1)^{\mu_{r-1}-\mu_r},\ldots,1^{\mu_1-\mu_2},0^{k+\mu_N-\mu_1},
(-1)^{\mu_{N-1}-\mu_N},\ldots,(r-N)^{-\mu_{r+1}})
\end{equation}

Furthermore, we write down the weight $\bw'=\sum_{i=1}^kw_i\omega_i-
\sum_{i=1}^kv_i\alpha_i$ as a sequence of integers
$(\sigma_1,\ldots,\sigma_k)$. The condition
$\fM_{(\zeta^\circ_\BC,\zeta^\bullet_\BR)}^{\on{reg}}(V,W)\ne\emptyset$
implies $\sigma_i\geq\sigma_{i+1}$ for $i\in I^0_0$, and
$\sigma_{k_0+\ldots+k_{p-1}+1}-\sigma_{k_0+\ldots+k_p}\leq N$
for any $0<p\leq n$,
where we put for convenience $k_0=0$.
The condition~\refe{uslovie} implies that
$\sigma_{k_1+\ldots+k_{p-1}+1}+\ldots+\sigma_{k_1+\ldots+k_p}=0$
for any $0<p\leq n$. Thus the sequence
$(\sigma_{k_1+\ldots+k_{p-1}+1},\ldots,\sigma_{k_1+\ldots+k_p})$
is a generalized
Young diagram to be denoted by $\sigma^{(p)}$. The transposed generalized Young
diagram $\blambda{}^{(p)}:={}^t\sigma^{(p)}$
corresponds to the same named level
$k_p$ dominant $\widehat{\fsl(N)}$-weight
$\blambda{}^{(p)}\in\Lambda^+_{k_p}(\widehat{\fsl(N)})$.

Recall that the affine Weyl group $W_{\on{aff}}$ acts on the set of level $N$
weights of $\widehat{\fsl(k)}$. If we write down these weights as the
sequences $(\chi_1,\ldots,\chi_k)$ then the action of $W_{\on{aff}}$ is
generated by permutations of $\chi_i$'s and the operations which only change
$\chi_i,\chi_j$ for some pair $i,j\in I$; namely, $\chi_i\mapsto\chi_i+N,\
\chi_j\mapsto\chi_j-N$.

\lem{otvr}
The sequence $(\sigma_1,\ldots,\sigma_k)$ is $W_{\on{aff}}$-conjugate to
$^t(\blambda{}^{(1)}+\ldots+\blambda{}^{(n)})$.
\elem

\prf
To simplify the notation we assume that $n=2$; the general case is not much
different.
Let $\blambda{}^{(1)}=(\lambda^{(1)}_1\geq\ldots\geq\lambda^{(1)}_N)$, and
$\blambda{}^{(2)}=(\lambda^{(2)}_1\geq\ldots\geq\lambda^{(2)}_N)$.
We set $\blambda=(\lambda_1\geq\ldots\geq\lambda_N)$ where
$\lambda_i:=\lambda^{(1)}_i+\lambda^{(2)}_i$. We assume
$\lambda^{(1)}_{r_1}>0\geq\lambda^{(1)}_{r_1+1}$ for some $0<r_1<N$, and
$\lambda^{(2)}_{r_2}>0\geq\lambda^{(2)}_{r_2+1}$ for some $0<r_2<N$.
If $r_1=r_2$, then the formula~\refe{javno} makes it clear that the sequence
$(\sigma_1,\ldots,\sigma_k)$ being a concatenation of the sequences
$(\sigma_1,\ldots,\sigma_{k_1})={}^t\blambda{}^{(1)}$ and
$(\sigma_{k_1+1},\ldots,\sigma_k)={}^t\blambda{}^{(2)}$ differs by a
permutation from the sequence $^t(\blambda{}^{(1)}+\blambda{}^{(2)})$.

Otherwise we assume $r_1>r_2$, and $\lambda_r>0\geq\lambda_{r+1}$ for some
$r_1\geq r\geq r_2$. Once again, to simplify the exposition, let us assume
that $r_1>r>r_2$.
According to the formula~\refe{javno}, if we reorder
the concatenation of $^t\blambda{}^{(1)}$ and $^t\blambda{}^{(2)}$ to obtain
a nonincreasing sequence, we get
$$(r_1^{\lambda^{(1)}_1},\ldots,r^{\lambda^{(1)}_r-\lambda^{(1)}_{r+1}},\ldots,
r_2^{\lambda^{(1)}_{r_2}-\lambda^{(1)}_{r_2+1}+\lambda^{(2)}_{r_2}},\ldots,
(r_1-N)^{-\lambda^{(1)}_{r_1+1}+\lambda^{(2)}_{r_1}-\lambda^{(2)}_{r_1+1}},
\ldots,$$ $$(r-N)^{\lambda^{(2)}_r-\lambda^{(2)}_{r+1}},\ldots,
(r_2-N)^{-\lambda^{(2)}_{r_2+1}})$$
On the other hand, the sequence $^t(\blambda{}^{(1)}+\blambda{}^{(2)})$ reads
$$
\begin{aligned}
(r^{\lambda^{(1)}_r+\lambda^{(2)}_r},\ldots,&
r_2^{\lambda^{(1)}_{r_2}+\lambda^{(2)}_{r_2}-\lambda^{(1)}_{r_2+1}-
\lambda^{(2)}_{r_2+1}},\ldots,\\
&(r_1-N)^{\lambda^{(1)}_{r_1}+\lambda^{(2)}_{r_1}-\lambda^{(1)}_{r_1+1}-
\lambda^{(2)}_{r_1+1}},\ldots,(r-N)^{-\lambda^{(1)}_{r+1}-\lambda^{(2)}_{r+1}})
\end{aligned}
$$
Now it is immediate to check that for any residue $h$ modulo $N$ its
multiplicity in the latter sequence is the sum of multiplicities of the same
residues in the former sequence. This means that the former sequence is
$W_{\on{aff}}$-conjugate to the latter one. The lemma is proved.
\epr

\ssec{dvo}{Birational convolution morphism}
Recall that we have a proper morphism $\pi_{0,\zeta^\bullet}:\
\fM_{(\zeta^\circ_\BC,\zeta^\bullet_\BR)}(V,W)
\to\fM_{(\zeta^\circ_\BC,0)}(V,W)$ introduced in \cite{Na}~3.2. Since $\bw'$ is
not necessarily dominant weight of $\widehat{\fsl(k)}$, the open
stratum $\fM_{(\zeta^\circ_\BC,0)}^{\on{reg}}(V,W)$ may be empty.
However, replacing $\bv$ by $\bv'=(v'_1,\ldots,v'_k)$ so that
$\bw''=\sum_{i=1}^kw''_i\omega_i:=\sum_{i=1}^kw_i\omega_i-
\sum_{i=1}^kv'_i\alpha_i$ is {\em dominant} and
$W_{\on{aff}}$-conjugate to $\bw'$, we can identify
$\fM_{(\zeta^\circ_\BC,0)}(V,W)$ with
$\fM_{(\zeta^\circ_\BC,0)}(V',W)$. Moreover, in this case the open subset
$\fM_{(\zeta^\circ_\BC,0)}^{\on{reg}}(V',W)$ is not empty, and the
morphism $\pi_{0,\zeta^\bullet}:\
\fM_{(\zeta^\circ_\BC,\zeta^\bullet_\BR)}(V,W)
\to\fM_{(\zeta^\circ_\BC,0)}(V',W)$ is birational. Recall that for
$\zeta^\circ_\BC=0$, in section~7 of~\cite{BF} we identified
$\fM_{(0,0)}(V',W)$ with the Uhlenbeck space
$\CU^\lambda_{\SL(N),\mu}(\BA^2/\Gamma_k)$ for certain level $k$
dominant $\fsl(N)_{\on{aff}}$-weights $\lambda,\mu$. In the
notations of current~\refss{transpon} we have $\mu=
(k,\bmu,-\frac{1}{2k}(2d+(\bmu,\bmu)-(\blambda,\blambda))),\
\lambda=(k,\blambda,0)$. Here $d=\sum_{i=1}^kv'_i$, and
$\blambda=\sum_{p=1}^n\blambda^{(p)}$ according to~\refl{otvr}.

For $1\leq p\leq n$ we introduce a level $k_p$ dominant
$\fsl(N)_{\on{aff}}$-weight $\lambda^{(p)}:=(k_p,\blambda^{(p)},0)$.
We set $\Blambda=(\lambda^{(1)},\ldots,\lambda^{(n)})$.
For arbitrary $\zeta^\circ_\BC$ we {\em define}
$\CU^\Blambda_{\SL(N),\mu}
(S_{\vec k}^{\vec a})$ as
$\fM_{(\zeta^\circ_\BC,\zeta^\bullet_\BR)}(V,W)$.
We {\em define} the convolution morphism
$\varpi:\ \CU_{\SL(N),\mu}^\Blambda
(S^{\vec a}_{\vec k})\to
\CU_{\SL(N),\mu}^\lambda(S'{}^{\vec a}_{\vec k})$ as
$\pi_{0,\zeta^\bullet}:\
\fM_{(\zeta^\circ_\BC,\zeta^\bullet_\BR)}(V,W)
\to\fM_{(\zeta^\circ_\BC,0)}(V',W)$. We will mostly use the particular case
$\varpi:\
\CU_{\SL(N),\mu}^\Blambda(S_{\vec k})\to
\CU_{\SL(N),\mu}^\lambda(S_k)=\CU_{\SL(N),\mu}^\lambda(\BA^2/\Gamma_k)$ defined as
$\pi_{0,\zeta^\bullet}:\
\fM_{(\zeta^\circ_\BC,\zeta^\bullet_\BR)}(V,W)
\to\fM_{(\zeta^\circ_\BC,0)}(V',W)$ for $\zeta^\circ_\BC=(0,\ldots,0)$.

\ssec{Tannaka}{Tensor product} Recall the notations
of~\refss{lambda}. Now the construction of~section~5(i)
of~\cite{NaM} gives rise to a morphism $\eta^\pm$ from
$\fM_{(\zeta^\circ_\BC,\zeta^\pm_\BR)}(V,W)$ to the moduli stack of
certain perverse coherent sheaves on
$\ol{S}{}^{\vec a}_{\vec k}$ trivialized at
$\ell_\infty$. It follows from~\refp{41}(1) (``only if'' part) that
the image of $\eta^-$ consists of torsion free sheaves, which
implies that the image of $\eta^+$ consists of the perverse sheaves
which are Serre-dual to the torsion free sheaves. We will denote the
connected component of the moduli stack of torsion free sheaves
(resp. of Serre-dual of torsion free sheaves) on
$\ol{S}{}^{\vec a}_{\vec k}$ birationally mapping to
$\CU_{\SL(N),\mu}^\Blambda
(S_{\vec k}^{\vec a})$ by $\CG
ies_\mu^\Blambda
(S_{\vec k}^{\vec a})$ (resp. by $\CS\CG
ies_\mu^\Blambda(S_{\vec k}^{\vec a})$).

\lem{Gies}
The morphisms $\eta^-:\ \fM_{(\zeta^\circ_\BC,\zeta^-_\BR)}(V,W)\to
\CG ies_\mu^\Blambda(S_{\vec k}^{\vec a}),\
\eta^+:\ \fM_{(\zeta^\circ_\BC,\zeta^+_\BR)}(V,W)\to
\CS\CG ies_\mu^\Blambda
(S_{\vec k}^{\vec a})$
are isomorphisms.
\elem

\prf
Follows from~\refp{41}(1) by the argument of~section~5 of~\cite{NaM}.
\epr

We consider the locally closed subvariety
$\bM_{(\zeta^\circ_\BC,\zeta^\bullet_\BR)}(V,W)\subset
\mu^{-1}(\zeta^\circ_\BC)\subset\bM(V,W)$ formed by all the
$\zeta^\bullet_\BR$-semistable modules. Let us denote by
$\Perv_\mu^\Blambda(S^{\vec a}_{\vec k})$
the moduli stack of perverse coherent sheaves on
$\ol{S}{}^{\vec a}_{\vec k}$ trivialized at $\ell_\infty$
and having the same numerical invariants as the torsion free sheaves in
$\CG ies_\mu^\Blambda
(S^{\vec a}_{\vec k})$.
The construction of~section~5(i) of~\cite{NaM}
gives rise to a morphism $\eta^\bullet$ from the stack
$\bM_{(\zeta^\circ_\BC,\zeta^\bullet_\BR)}(V,W)/GL_V$ to
$\Perv_\mu^\Blambda(S^{\vec a}_{\vec k})$.

\lem{per}
$\eta^\bullet:\ \bM_{(\zeta^\circ_\BC,\zeta^\bullet_\BR)}(V,W)/GL_V\to
\Perv_\mu^\Blambda(S^{\vec a}_{\vec k})$
is an isomorphism.
\elem

\prf
Follows from~\refp{41}(2) by the argument of~section~5 of~\cite{NaM}.
\epr

It follows from~\refl{Gies} that
we have a projective morphism $\pi_{\zeta^\bullet,\zeta^-}:\
\CG ies_\mu^\Blambda(S^{\vec a}_{\vec k})
\to\CU_{\SL(N),\mu}^\Blambda
(S^{\vec a}_{\vec k})$
(resp. $\pi_{\zeta^\bullet,\zeta^+}:\
\CS\CG ies_\mu^\Blambda
(S^{\vec a}_{\vec k})\to
\CU_{\SL(N),\mu}^\Blambda
(S^{\vec a}_{\vec k})$).

\lem{perv}
If $E$ is a torsion free coherent sheaf on
$S^{\vec a}_{\vec k}$,
and $E'$ is the Serre dual of a torsion free coherent sheaf on
$S^{\vec a}_{\vec k}$,
then $E\otimes E'$ is a perverse coherent sheaf on
$S^{\vec a}_{\vec k}$.
\elem

\prf
Clearly, $\unl{H}^{>1}(E\otimes E')$ vanishes,
$\unl{H}^1(E\otimes E')$ is a torsion sheaf supported at finitely many points,
and $\unl{H}^0(E\otimes E')$ is torsion free. The same is true for the Serre
dual sheaf of $E\otimes E'$ (being a tensor product of the same type).
\epr

Thus we obtain a morphism $\CG
ies_\mu^\Blambda(S^{\vec a}_{\vec k})
\times\CS\CG ies_{\mu'}^{{}'\Blambda}
(S^{\vec a}_{\vec k})\to
\Perv_{\mu\otimes\mu'}^{\Blambda\otimes{}'\Blambda}(S^{\vec a}_{\vec k})$.
Here we understand $\bmu$ (resp. $\bmu{}'$) as a homomorphism
$\Gamma_k\to\SL(N)$ (resp. $\Gamma_k\to\SL(N')$), and
$\bmu\otimes\bmu{}'$ as the tensor product homomorphism
$\Gamma_k\to\SL(NN')$; similarly for $\blambda$'s. Furthermore, we
set $\lambda^{(p)}\otimes{}'\lambda^{(p)}:=
(k_i,\blambda{}^{(p)}\otimes{}'\blambda{}^{(p)},0)$, and
$\Blambda\otimes{}'\Blambda=(\lambda^{(1)}\otimes{}'\lambda^{(1)},\ldots,
\lambda^{(n)}\otimes{}'\lambda^{(n)})$.
Finally, for $\mu=(k,\bmu,m),\ \mu'=(k,\bmu{}',m')$ we set
$\mu\otimes\mu':=(k,\bmu\otimes\bmu{}',{\mathsf m})$ where
$${\mathsf m}:=mN'+m'N+\frac{1}{2k}\left[N'(\bmu,\bmu)-
N'(\sum_{p=1}^n\blambda{}^{(p)},\sum_{p=1}^n\blambda{}^{(p)})+
N(\bmu{}',\bmu{}')-\right.$$
$$\left.-N(\sum_{p=1}^n{}'\blambda{}^{(p)},\sum_{p=1}^n{}'\blambda{}^{(p)})-
(\bmu\otimes\bmu{}',\bmu\otimes\bmu{}')+
(\sum_{p=1}^n\blambda{}^{(p)}\otimes{}'\blambda{}^{(p)},
\sum_{p=1}^n\blambda{}^{(p)}\otimes{}'\blambda{}^{(p)})\right].$$

Composing this morphism with the further
projection (due to~\refl{per})
$\Perv_{\mu\otimes\mu'}^{\Blambda\otimes{}'\Blambda}(S^{\vec a}_{\vec k})\to
\CU_{\SL(NN'),\mu\otimes\mu'}^{\Blambda\otimes{}'\Blambda}(S^{\vec a}_{\vec k})$
we obtain the morphism
$\tau:\ \CG ies_\mu^\Blambda
(S^{\vec a}_{\vec k})\times
\CS\CG ies_{\mu'}^{{}'\Blambda}
(S^{\vec a}_{\vec k})\to
\CU_{\SL(NN'),\mu\otimes\mu'}^{\Blambda\otimes{}'\Blambda}(S^{\vec a}_{\vec k})$.


\prop{tens}
The morphism $\tau$ factors through
$\bar\tau:\
\CU_{\SL(N),\mu}^\Blambda
(S^{\vec a}_{\vec k})\times
\CU_{\SL(N'),\mu'}^{{}'\Blambda}
(S^{\vec a}_{\vec k})\to
\CU_{\SL(NN'),\mu\otimes\mu'}^{\Blambda\otimes{}'\Blambda}(S^{\vec a}_{\vec k})$.
\eprop

\prf
Let us denote $\CG ies_\mu^\Blambda
(S^{\vec a}_{\vec k})\times
\CS\CG ies_{\mu'}^{{}'\Blambda}
(S^{\vec a}_{\vec k})$ by $X$, and
$\CU_{\SL(N),\mu}^\Blambda
(S^{\vec a}_{\vec k})\times
\CU_{\SL(N'),\mu'}^{{}'\Blambda}
(S^{\vec a}_{\vec k})$ by $Y$,
and $\CU_{\SL(NN'),\mu\otimes\mu'}^{\Blambda\otimes{}'\Blambda}(S^{\vec a}_{\vec k})$ by $Z$
for short. We have to prove that the morphism $\tau:\ X\to Z$ factors through
the morphism $\pi:=\pi_{\zeta^\bullet_\BR,\zeta^-_\BR}\times
\pi_{\zeta^\bullet_\BR,\zeta^+_\BR}:\ X\to Y$, and a morphism $\bar\tau:\
Y\to Z$. It is easy to see that $\tau$ contracts the fibers of $\pi$, that is
for any $y\in Y$ we have $\tau(\pi^{-1}(y))=z$ for a certain point $z\in Z$.
It means that the image $T$ of $\pi\times\tau:\ X\to Y\times Z$ projects onto
$Y$ bijectively. Furthermore, $T$ is a closed
subvariety of $Y\times Z$ since both $\pi$ and $\tau$ are proper. Finally,
$Y$ is normal by a theorem of Crawley-Boevey. This implies that the projection
$T\to Y$ is an isomorphism of algebraic varieties.
Hence $T$ is the graph of a
morphism $Y\to Z$. This is the desired morphism $\bar\tau$.

This argument was explained to us by A.~Kuznetsov.
\epr

\sec{cotan}{Tannakian approach}

\ssec{tan}{} Given an almost simple simply connected group $G$, and the weights
$\bmu\in\Lambda^+_k,\ \blambda{}^{(i)}\in\Lambda^+_{k_i},\ 1\leq i\leq n$, and a
positive integer $d$, we consider the moduli space
$\Bun_{G,\bmu}^{\Bblambda,d/k}
(S^{\vec a}_{\vec k})$ introduced in~\refss{BMKS}. It classifies
the $G$-bundles on the stack $\ol{S}{}^{\vec a}_{\vec k}$ of
second Chern class $d/k$,
trivialized at infinity such that the class of the fiber at infinity is given
by $\bmu$, while the class of the fiber at $s_i$ is given by $\blambda^{(i)}$.

\conj{connect}
$\Bun_{G,\bmu}^{\Bblambda,d/k}(S_{\vec k})$
is connected (possibly empty). \econj

Following the numerology of~\refss{dvo}, we introduce the weights
$\lambda^{(i)}:=(k_i,\blambda{}^{(i)},0)\in\Lambda^+_{\on{aff},k_i}$, and
$\mu:=(k,\bmu,-\frac{1}{2k}(2d+(\bmu,\bmu)-(\blambda,\blambda)))\in
\Lambda^+_{\on{aff},k}$ where $\blambda=\sum_{i=1}^n\blambda{}^{(i)}$.
We also set $\lambda:=(k,\blambda,0)$. Now we redenote
$\Bun_{G,\bmu}^{\Bblambda,d/k}
(S^{\vec a}_{\vec k})$ by
$\Bun_{G,\mu}^\Blambda(S^{\vec a}_{\vec k})$.

Given a representation $\varrho:\ G\to\SL(W_\varrho)$ we have a morphism
$\varrho_*:\ \Bun_{G,\mu}^\Blambda
(S^{\vec a}_{\vec k})\to
\Bun_{\SL(W_\varrho),\varrho\circ\mu}^{\varrho\circ\Blambda}(S^{\vec a}_{\vec k})
\subset
\CU_{\SL(W_\varrho),\varrho\circ\mu}^{\varrho\circ\Blambda}(S^{\vec a}_{\vec k})$.
Here $\lambda^{(p)}=(k_p,\blambda{}^{(p)},0),\
\varrho\circ\lambda^{(p)}:=(k_p,\varrho\circ\blambda{}^{(p)},0);\\
\varrho\circ\Blambda=
(\varrho\circ\lambda^{(1)},\ldots,\varrho\circ\lambda^{(n)});\
\mu=(k,\bmu,m),\
\varrho\circ\mu:=(k,\varrho\circ\bmu,\varrho_\BZ m)$, and $\varrho_\BZ$ is
the Dynkin index of the representation $\varrho$ (we stick to the notation
of~\cite{BFG},~6.1).

We define $\CU_{G,\mu}^\Blambda
(S^{\vec a}_{\vec k})$ as the closure of the image
of
$\prod_\varrho\varrho_*(\Bun_{G,\mu}^\Blambda
(S^{\vec a}_{\vec k}))$ inside
$\prod_\varrho\CU_{\SL(W_\varrho),\varrho\circ\mu}
^{\varrho\circ\Blambda}(S^{\vec a}_{\vec k})$.
For any $\varrho$ we have an evident projection morphism
$\CU_{G,\mu}^\Blambda
(S^{\vec a}_{\vec k})\to\CU_{\SL(W_\varrho),\varrho\circ\mu}
^{\varrho\circ\Blambda}(S^{\vec a}_{\vec k})$. By
an abuse of notation we will denote this morphism $\varrho_*$.

\prop{closed} Assume that any representation of $G$ is a direct
summand of a tensor power of $\varrho$ (this is equivalent to requesting that
$\varrho$ is faithful). Then $\varrho_*:\
\CU_{G,\mu}^\Blambda
(S^{\vec a}_{\vec k})\to\CU_{\SL(W_\varrho),\varrho\circ\mu}
^{\varrho\circ\Blambda}(S^{\vec a}_{\vec k})$ is
a closed embedding. In particular,
$\CU_{G,\mu}^\Blambda
(S^{\vec a}_{\vec k})$ is of finite type. \eprop

\prf Let $x\in\CU_{\SL(W_\varrho),\varrho\circ\mu}
^{\varrho\circ\Blambda}(S^{\vec a}_{\vec k})$ be
a point in the closure of the locally closed subvariety
$\varrho_*(\Bun_{G,\mu}^\Blambda
(S^{\vec a}_{\vec k}))$. There is an affine pointed
curve $(C,c)\subset\CU_{\SL(W_\varrho),\varrho\circ\mu}
^{\varrho\circ\Blambda}(S^{\vec a}_{\vec k})$
such that
$(C-c)\subset\varrho_*(\Bun_{G,\mu}^\Blambda
(S^{\vec a}_{\vec k}))$, and $c=x$.

Let $\varsigma:\ G\to\SL(W_\varsigma)$ be another representation of
$G$. We choose a projection $\varkappa:\ \varrho^{\otimes
m}\twoheadrightarrow\varsigma$. According to~\refp{tens}, we
consider ${\bar\tau}(C)\subset \CU_{\SL(W_\varrho^{\otimes
m}),\varrho^{\otimes m}\circ\mu} ^{\varrho^{\otimes
m}\circ\Blambda}
(S^{\vec a}_{\vec k})$, and then
$\varkappa_*{\bar\tau}(C)\subset
\CU_{\SL(W_\varsigma),\varsigma\circ\mu}
^{\varsigma\circ\Blambda}(S^{\vec a}_{\vec k})$.
Since $\varkappa_*{\bar\tau}(C-c)\subset
\Bun_{\SL(W_\varsigma),\varsigma\circ\mu}
^{\varsigma\circ\Blambda}
(S^{\vec a}_{\vec k})$, we have lifted $x=c$ to a
point of $\CU_{G,\mu}^\Blambda
(S^{\vec a}_{\vec k})$.

It remains to prove that such a lift is unique. Let $\varkappa':\
\varrho^{\otimes m'}\twoheadrightarrow\varsigma$ be another
projection. Then $\varkappa_*{\bar\tau}=\varkappa'_*{\bar\tau}:\
(C-c)\hookrightarrow\Bun_{\SL(W_\varsigma),\varsigma\circ\mu}
^{\varsigma\circ\Blambda}(S^{\vec a}_{\vec k})$.
Since $\CU_{\SL(W_\varsigma),\varsigma\circ\mu}
^{\varsigma\circ\Blambda}(S^{\vec a}_{\vec k})$
is separated, it follows that
$\varkappa_*{\bar\tau}(c)=\varkappa'_*{\bar\tau}(c)$.

This completes the proof of the proposition.
\epr

\ssec{convo}{Convolution morphism}
The collection of convolution morphisms  $\CU_{\SL(W_\varsigma),\varsigma\circ\mu}
^{\varsigma\circ\Blambda}(S^{\vec a}_{\vec k})\to
\CU_{\SL(W_\varsigma),\varsigma\circ\mu}^{\varsigma\circ\lambda^{(1)}+\ldots+
\varsigma\circ\lambda^{(n)}}(S'{}^{\vec a}_{\vec k})$ (see~\refss{dvo}) gives rise to
the convolution morphism $\varpi:\ \CU_{G,\mu}^\Blambda
(S^{\vec a}_{\vec k})\to\CU_{G,\mu}^{\lambda^{(1)}+\ldots+\lambda^{(n)}}
(S'{}^{\vec a}_{\vec k})$.

\lem{birat} The morphism $\varpi:\ \CU_{G,\mu}^\Blambda
(S^{\vec a}_{\vec k})\to\CU_{G,\mu}^{\lambda^{(1)}+\ldots+\lambda^{(n)}}
(S'{}^{\vec a}_{\vec k})$ is birational.
\elem

\prf
It suffices to check that $\varpi$ is an isomorphism when restricted to
the open subset $\Bun_{G,\mu}^{\lambda^{(1)}+\ldots+\lambda^{(n)}}
(S'{}^{\vec a}_{\vec k})\subset\CU_{G,\mu}^{\lambda^{(1)}+\ldots+\lambda^{(n)}}
(S'{}^{\vec a}_{\vec k})$. For any representation $\varsigma:\ G\to\SL(W_\varsigma)$
and the corresponding closed embedding $\varsigma_*:\
\CU_{G,\mu}^{\lambda^{(1)}+\ldots+\lambda^{(n)}}(S'{}^{\vec a}_{\vec k})\hookrightarrow
\CU_{\SL(W_\varsigma),\varsigma\circ\mu}^{\varsigma\circ\lambda^{(1)}+
\ldots+\varsigma\circ\lambda^{(n)}}(S'{}^{\vec a}_{\vec k})$ we have
$\varsigma_*(\Bun_{G,\mu}^{\lambda^{(1)}+\ldots+\lambda^{(n)}}(S'{}^{\vec a}_{\vec k}))\subset
\Bun_{\SL(W_\varsigma),\varsigma\circ\mu}^{\varsigma\circ\lambda^{(1)}+
\ldots+\varsigma\circ\lambda^{(n)}}(S'{}^{\vec a}_{\vec k})$. Any
vector bundle
$\CF\in\Bun_{\SL(W_\varsigma),\varsigma\circ\mu}^{\varsigma\circ\lambda^{(1)}+
\ldots+\varsigma\circ\lambda^{(n)}}(S'{}^{\vec a}_{\vec k})$
 has a unique preimage
$\CF'$ under the convolution morphism $\CU_{\SL(W_\varsigma),\varsigma\circ\mu}
^{\varsigma\circ\Blambda}(S^{\vec a}_{\vec k})\to
\CU_{\SL(W_\varsigma),\varsigma\circ\mu}^{\varsigma\circ\lambda^{(1)}+\ldots+
\varsigma\circ\lambda^{(n)}}(S'{}^{\vec a}_{\vec k})$, moreover,
$\CF'\in\Bun_{\SL(W_\varsigma),\varsigma\circ\mu}
^{\varsigma\circ\Blambda}(S^{\vec a}_{\vec k})$ since the
0-stability implies the $\zeta^\bullet_\BR$-stability. Clearly, if
$\CF$ lies in the image
$\varsigma_*(\Bun_{G,\mu}^{\lambda^{(1)}+\ldots+\lambda^{(n)}}(S'{}^{\vec a}_{\vec k}))\subset
\Bun_{\SL(W_\varsigma),\varsigma\circ\mu}^{\varsigma\circ\lambda^{(1)}+
\ldots+\varsigma\circ\lambda^{(n)}}(S'{}^{\vec a}_{\vec k})$, then $\CF'$ lies in the
image $\varsigma_*(\Bun_{G,\mu}^\Blambda
(S^{\vec a}_{\vec k}))\subset\Bun_{\SL(W_\varsigma),\varsigma\circ\mu}
^{\varsigma\circ\Blambda}(S^{\vec a}_{\vec k})$. So this
$\CF'$ is the unique preimage of $\CF$ under $\varpi$.
\epr

\ssec{main}{Main conjecture}
We have a proper surjective morphism
$\varpi:\ \CU_{G,\mu}^\Blambda(S_{\vec k})\to
\CU_{G,\mu}^{\lambda^{(1)}+\ldots+\lambda^{(n)}}(\BA^2/\Gamma_k)$, and we are
interested in the multiplicities in
$\varpi_*IC(\CU_{G,\mu}^\Blambda(S_{\vec k}))$.
For $\mu\leq\nu\leq\lambda^{(1)}+\ldots+\lambda^{(n)}$ we
will denote the multiplicity of $IC(\CU_{G,\mu}^\nu(\BA^2/\Gamma_k))$ in
$\varpi_*IC(\CU_{G,\mu}^\Blambda(S_{\vec k}))$
by $M_{\nu,\mu}^\Blambda$.

\conj{key}
a) $\varpi$ is semismall, and hence
$M_{\nu,\mu}^\Blambda$ is just a vector space in
degree zero.

b) $M_{\nu,\mu}^\Blambda$ is independent of $\mu$ and
equals the multiplicity of the $G_{\aff}^{\vee}$-module $L(\nu)$ in
the tensor product $L(\lambda^{(1)})\otimes\ldots\otimes L(\lambda^{(n)})$.
\econj

\rem{hiraku}
The direct image
$\varpi_*IC(\CU_{G,\mu}^\Blambda(S_{\vec k}))$
{\em is not} isomorphic to the direct sum
$\bigoplus_\nu M_{\nu,\mu}^\Blambda\otimes
IC(\CU_{G,\mu}^\nu(\BA^2/\Gamma_k))$: it contains the IC sheaves of other
strata of $\CU_{G,\mu}^\nu(\BA^2/\Gamma_k)$ with nonzero multiplicities.
This observation is due to H.~Nakajima~\cite{Na}.
\erem

The following Proposition is essentially proved in~\cite{Na}.

\prop{nakaji}
~\refco{key} for $G=\SL(N)$ holds true.
\eprop

\prf By definition, the desired multiplicity
$M_{\nu,\mu}^\Blambda$ can be computed on the
quiver varieties, in the particular case
$\zeta^\circ_\BC=(0,\ldots,0)$. It is computed in Theorem~5.15
(equation~(5.16)) of~\cite{Na} under the name $V^{{\mathbf
v}^0,\emptyset}_{{\mathbf v}',\emptyset}$. Note that we are
interested in the particular case $\mu=\emptyset=\lambda,\ {\mathbf
v}^0=\bv,$ thereof (we apologize for the conflicting roles of
$\lambda,\mu$ in {\em loc. cit.} and in the present paper). We set
$d=v_1+\ldots+v_k,\ d'=v'_1+\ldots+v'_k$. Finally, $\bnu$ is
associated to the pair $({\mathbf v}',{\mathbf w})$ as
in~\cite{BF}~7.3, and
$\nu=\left(k,\bnu,\frac{1}{2k}\left[2d'-2d-(\bnu,\bnu)+
(\blambda^{(1)}+\ldots+\blambda^{(n)},
\blambda^{(1)}+\ldots+\blambda^{(n)})\right]\right)$.

Furthermore, in the~Remark~5.17.(3) of~\cite{Na} the multiplicity
$V^{{\mathbf v}^0,\emptyset}_{{\mathbf v}',\emptyset}$ is identified
via I.~Frenkel's level-rank duality with the multiplicity of the
$G_{\aff}^{\vee}$-module $L(\nu)$ in the tensor product
$L(\lambda^{(1)})\otimes\ldots\otimes L(\lambda^{(n)})$. \epr
\ssec{fam-finite}{Digression on the Beilinson-Drinfeld Grassmannian}
Let $C$ be a smooth algebraic curve and let $c$ be a point of
$C$. It is well-known that a choice of formal parameter at $c$ gives rise to an identification
of $\Gr_G$ with the moduli space of $G$-bundles on $C$ endowed with a trivialization away from
$c$. Similarly, for any $n\geq 1$ one can introduce the {\em Beilinson-Drinfeld Grassmannian}
$\Gr_{C,G,n}$ as the moduli space of the following data:

1) An ordered collection of points $(c_1,\ldots,c_n)\in C^n$;

2) A $G$-bundle $\calF$ on $C$ trivialized away from $(c_1,\ldots,c_n)$.

We have an obvious map $p_n:\Gr_{C,G,n}\to C^n$ sending the above data
to $(c_1,\ldots,c_n)$. When all the points $c_i$ are distinct, the fiber
$p_n^{-1}(c_1,\ldots,c_n)$ is non-canonically isomorphic to $(\Gr_G)^n$. When
all the points coincide, the corresponding fiber is isomorphic
to just one copy of $\Gr_G^n$. For any $\lam_1,\ldots,\lam_n\in\Lam^+$
one can define the closed subvariety $\oGr_{C,G}^{\lam_1,\ldots,\lam_n}$
in $\Gr_{C,G,n}$ such that for any collection $(p_1,\ldots,p_n)$ of distinct
points of $C$ the intersection $p_n^{-1}(c_1,\ldots,c_n)\cap
\oGr_{C,G}^{\lam_1,\ldots,\lam_n}$
is isomorphic to $\oGr^{\lam_1}\x\ldots\x\oGr^{\lam_n}$ and the intersection
$p_n^{-1}(c,\ldots,c)\cap \oGr_{C,G}^{\lam_1,\ldots,\lam_n}$ is isomorphic
to $\oGr^{\lam_1+\ldots+\lam_n}$.

Similarly, given $C$ and $n$ as above one defines
the scheme $\tGr_{G,C,n}$ classifying the following
data:

1) An element $(c_1,\ldots,c_n)\in C^n$.

2) An $n$-tuple $(\calF_1,\ldots,\calF_n)$ of $G$-bundles on $C$; we also let
$\calF_0$ denote the trivial $G$-bundle on $C$.

3) An isomorphism $\kap_i$ between $\calF_{i-1}|_{C\backslash \{ c_i\}}$ and $\calF_{i}|_{C\backslash \{ c_i\}}$
for each $i=1,\ldots,n$.

\noindent
We denote by $\tilp_n$ the natural map from $\tGr_{C,G,n}$ to $C^n$.
Note that from 3) one gets a trivialization of $\calF_n$ away from $(c_1,\ldots,c_n)$.
Thus we have the natural map $\tGr_{C,G,n}\to \Gr_{C,G,n}$. This map is proper and it is an isomorphism
on the open subset where all the points $c_i$ are distinct.
On the other hand, the morphism $\tilp_n^{-1}(c,\ldots,c)\to \Gr_{C,G,n}$ is isomorphic to the morphism
$\underset{n}{\underbrace{\Gr_G\star\cdots\star\Gr_G}}\to \Gr_G$.
For $\lam_1,\ldots,\lam_n$ as above we denote by $\tGr_{C,G}^{\lam_1,\ldots,\lam_n}$
the closed subset of $\tGr_{C,G,n}$ given by the condition that each $\kap_i$ lies $\oGr_G^{\lam_i}$.
Then the intersection $\tilp_n^{-1}(c,\ldots,c)\cap \tGr_{C,G}^{\lam_1,\ldots,\lam_n}$ is
isomorphic to $\oGr_G^{\lam_1}\star\ldots\star\oGr_G^{\lam_n}$.
\ssec{fam}{Beilinson-Drinfeld Grassmannian for $G_{\aff}$} Our next task will
be to define an analog of (some pieces) of the Beilinson-Drinfeld Grassmannian
for $G_{\aff}$ in the case when $C=\AA^1$. The idea is that as $(a_1,\ldots,a_n)\in\BA^{n-1}$ varies, we will
organize $\CU_{G,\mu}^\Blambda
(S^{\vec a}_{\vec k})$ (resp.
$\CU_{G,\mu}^\Blambda
(S'{}^{\vec a}_{\vec k})$) into a family
$\CU_{G,\mu}^\Blambda(\CS_{\vec k})$
(resp.
$\CU_{G,\mu}^\Blambda(\CS'_{\vec k})$)
over $\CX=\BA^{n-1}$ (though there is {\em no} family of smooth
2-dimensional stacks over $\CX$). We will also construct a proper birational
morphism
$\varpi:\ \CU_{G,\mu}^\Blambda(\CS_{\vec k})\to
\CU_{G,\mu}^\Blambda
(\CS'_{\vec k})$ specializing to the morphisms $\varpi$
of~\refss{convo} for the particular values of $(a_1,\ldots,a_n)$.

In case $G=\SL(N)$, we {\em define}
$\CU_{G,\mu}^\Blambda(\CS_{\vec k})$
(resp.
$\CU_{G,\mu}^\Blambda(\CS'_{\vec k})$)
as the families of quiver varieties $\fN_{\zeta^\bullet_\BR}(V,W)$
(resp. $\fN_0(V,W)$) over the variety $\CX$ of moment levels
$\zeta^\circ_\BC$ (recall that $\zeta^\circ_\BC$ is reconstructed
from $(a_1,\ldots,a_n)$ by the beginning of~\refss{zeta}),
see~\cite{Na} between~Lemma~5.12 and~Remark~5.13.

For general $G$ we repeat the procedure of~\refss{tan}. We only have to
define the morphism $\bar\tau:\
\CU_{\SL(N),\mu}^\Blambda
(\CS^{\vec a}_{\vec k})\times
\CU_{\SL(N'),\mu'}^{{}'\Blambda}
(\CS^{\vec a}_{\vec k})\to
\CU_{\SL(NN'),\mu\otimes\mu'}^{\Blambda\otimes{}'\Blambda}(\CS^{\vec a}_{\vec k})$,
that is to prove a relative analogue of~\refp{tens}.

To this end we consider the resolution
$\fN_{\zeta_\BR}(V,W)\to\fN_{\zeta^\circ_\BR}(V,W)\to\fN_{\zeta^\bullet_\BR}(V,W)$,
see~\cite{Na} between~Lemma~5.12 and~Remark~5.13.
Here $\zeta^\circ_\BR$ is (chosen and fixed) generic in the
hyperplane $\zeta^\circ_\BR\cdot\delta=0$ (see~\cite{NaM}~1(iii)), and
$\zeta_\BR$ is in the chamber containing $\zeta^\circ_\BR$ in its closure
with $\zeta_\BR\cdot\delta<0$. According to the
Main~Theorem of~\cite{NaM}, $\fN_{\zeta_\BR}(V,W)$ is isomorphic to
the Giesecker moduli space of torsion-free sheaves on the
simultaneous resolution $\widetilde{\ol\CS}_{\vec k}$
trivialized at $\ell_\infty$. Now repeating the argument
of~\refp{tens} we obtain a morphism $\tilde\tau:\
\fN_{\zeta^\circ_\BR}(V,W)\times\fN_{\zeta^\circ_\BR}(V',W')\to
\fN_{\zeta^\circ_\BR}(V'',W'')$ where the $\Gamma_k$-modules $V'',W''$ are
defined as follows: $W''=W\otimes W',\ V''=V\otimes W'\oplus V'\otimes W\oplus
V\otimes V'\otimes(Q\ominus\BC^2)$, and $Q$ is the tautological 2-dimensional
representation of $\Gamma_k\subset\SL(2)$, while $\BC^2$ is the trivial
2-dimensional representation of $\Gamma_k$.
Composing it with the projection $\fN_{\zeta^\circ_\BR}(V'',W'')\to
\fN_{\zeta^\bullet_\BR}(V'',W'')$
we obtain a morphism $\tau':\
\fN_{\zeta^\circ_\BR}(V,W)\times\fN_{\zeta^\circ_\BR}(V',W')\to
\fN_{\zeta^\bullet_\BR}(V'',W'')$.
Now the argument of~\refp{tens} proves that $\tau'$ factors through
the desired morphism $\bar\tau:\
\fN_{\zeta^\bullet_\BR}(V,W)\times\fN_{\zeta^\bullet_\BR}(V',W')\to
\fN_{\zeta^\bullet_\BR}(V'',W'')$.


\ssec{taub-nut}{Semismallness of $\varpi$} In this section we speculate on a
possible approach to~\refco{key}a) using the double affine version of the
Beilinson-Drinfeld Grassmannian.

Assume $k$ is even. We
set $\CX=\BA^{n-1}$ with coordinates $a_1,\ldots,a_n,\
k_1a_1+\ldots+k_na_n=0$. We consider the weighted projective space
$\unl\BP(2,k,k,2)= (\BA^4-0)/\!/{\mathbb G}_m$ where ${\mathbb G}_m$
acts as follows: $t(x,y,z,w)=(t^2x,t^ky,t^kz,t^2w)$. We define a
relative surface $q:\ \unl{\ol\CS{}'}\to\CX$ as the hypersurface in
$\unl\BP(2,k,k,2)\times\CX$ given by the equation
$yz=(x-a_1w)^{k_1}\ldots(x-a_nw)^{k_n}$. The divisor at infinity is
given by $w=0$. Note that $q$ is a compactification of a subfamily
of the semiuniversal deformation of the $A_{k-1}$-singularity
constructed in~\cite{Sa}. Clearly, the fiber $q^{-1}(a_1,\ldots,a_n)$
with the divisor at infinity removed is isomorphic to
$\unl{S'}{}^{\vec a}_{\vec k}$. There is a family $\unl{\ol\CS}
\stackrel{p}{\to}\unl{\ol\CS{}'}\stackrel{q}{\to}\CX$ such that $p$ is
an isomorphism in a neighbourhood of the divisor at infinity, and the
restriction of $p$ to the fiber $q^{-1}(a_1,\ldots,a_n)$
with the divisor at infinity removed is nothing else than the partial
resolution $\unl{S}{}^{\vec a}_{\vec k}\to
\unl{S'}{}^{\vec a}_{\vec k}$ of~\refss{zeta}.

By the axioms of~\refss{axio}, we should have a proper birational
morphism $\varpi:\ \CU_G(\CS)\to\CU_G(\CS')$ whose fiber over $x=(0,\ldots,0)
\in\CX$ coincides with $\varpi$ of~\refss{convo}. This is nothing else than
$\varpi$ of~\refss{fam}.

Since the family $\unl{\ol\CS}\to\CX$ is equisingular, we expect the morphism
${\mathfrak p}:\ \CU_G(\CS)\to\CX$ to be
locally acyclic. Hence the specialization of the intersection cohomology sheaf
$\IC(\CU_G(\CS))$ to the fiber ${\mathfrak p}^{-1}(0,\ldots,0)$ coincides with
$\IC(\CU_G(S_{\vec k}))$. Since the specialization commutes with the
direct image under proper morphisms, we obtain
$\varpi_*\IC(\CU_G(S_{\vec k}))=
{\mathbf{Sp}}_{(0,\ldots,0)}\varpi_*\IC(\CU_G(\CS))=
{\mathbf{Sp}}_{(0,\ldots,0)}\IC(\CU_G(\CS'))$. Here the second equality holds
since $\varpi$ is an isomorphism off the diagonals in $\CX$. It follows that
$\varpi_*\IC(\CU_G(S_{\vec k}))$ is perverse (and semisimple, by the
decomposition theorem).

\bigskip
\footnotesize{
{\bf A.B.}: Department of Mathematics, Brown University,
151 Thayer St., Providence RI
02912, USA;\\
{\tt braval@math.brown.edu}}

\footnotesize{
{\bf M.F.}: IMU, IITP and National Research University Higher School of
Economics\\
Department of Mathematics, 20 Myasnitskaya st, Moscow 101000 Russia;\\
{\tt fnklberg@gmail.com}}

\end{document}